\documentclass[10pt,a4paper]{article}
\usepackage[left=1.00cm, right=1.00cm, top=2cm, bottom=2.0cm]{geometry}

\usepackage{graphicx,xcolor}
\usepackage{ifthen,calc}
\usepackage{amsfonts,amssymb,amsmath,amsthm}
\usepackage[colorlinks,breaklinks]{hyperref}
\hypersetup{urlcolor=blue, citecolor=red}

\usepackage{float}
\usepackage{indentfirst}

\usepackage{txfonts}

\usepackage{paralist}
\usepackage{authblk}

\newcommand{\review}[1]{\textcolor{black}{#1}}

\newcommand{\var}{\ensuremath{\mathrm{var}}}
\numberwithin{equation}{section}

\usepackage{units,comment, subcaption, color}
\usepackage[normalem]{ulem}

\newtheorem{remark}{Remark}

\usepackage[pagewise]{lineno}
\usepackage{paralist}

 \usepackage{etoolbox}

\begin{document}

\title{Uncertainty propagation and sensitivity analysis: results from the Ocular Mathematical Virtual Simulator}

\date{}

\author[$1$]{C. Prud'homme}
\author[$2,*$]{L. Sala}
\author[$3$]{M. Szopos}

\affil[$1$]{IRMA UMR CNRS 7501, Université de Strasbourg, Strasbourg 67000, France}
\affil[$2$]{Centre de recherche INRIA de Paris, Paris 75012, France}
\affil[$3$]{MAP5 UMR CNRS 8145, Universit\'e de Paris, Paris 75006, France}
\affil[$*$]{\textit{lorenzo.sala@inria.fr}}

\maketitle

\vspace{-3mm}
\begin{center}
\begin{minipage}{0.9\linewidth}		
\textbf{Abstract.}
We propose an uncertainty propagation study and a sensitivity analysis with the Ocular Mathematical Virtual Simulator, a computational and mathematical model that predicts the hemodynamics and biomechanics within the human eye.
In this contribution, we focus on the effect of intraocular pressure, retrolaminar tissue pressure and systemic blood pressure on the ocular posterior tissue vasculature.
The combination of a physically-based model with experiments-based stochastic input allows us to gain a better understanding of the physiological system, accounting both for the driving mechanisms and the data variability. \\

\textbf{Keywords:} predictive ocular vascular dynamics, computational model, sensitivity analysis, ocular mathematical virtual simulator, Sobol index analysis, uncertainty quantification
\end{minipage}
\end{center}

\section{Introduction}
\label{sec:intro}

The interest in patient-specific mathematical models applied to biomedical problems has greatly increased in the last years. 
In particular, the need for a better understanding and knowledge of quantities in the medical context has raised tremendously the complexity of the mathematical models employed to describe such physical systems. 
However, a crucial aspect to guarantee that the model and its numerical solutions are meaningful from the biomedical viewpoint is how inherent uncertainties are incorporated, {as recently discussed for instance in~\cite{hose2019cardiovascular}.}
In this direction, several works that studied the impact of uncertainties in the domain of cardiovascular
disease modelling showed particular promise for elucidating the complex interplay between hemodynamics, biomechanics, and electrophysiology. 
Examples include arterial hemodynamics~\cite{chen2013,Eck2015,Leguy2011,brault2017}, cardiovascular simulations~\cite{sankaran2011,marquis2018practical,fleeter2020multilevel}, electrophysiology~\cite{quaglino2018}, possibly coupled with electromechanical simulations~\cite{hurtado2017} and/or hemodynamics~\cite{campos2020}.\\

To the best of our knowledge, the eye's mathematical and computational modelling is still at its early stages, as recently reviewed in~\cite{guidoboni2019}. 
Biomechanical and fluid-dynamical aspects are of particular relevance for several clinical conditions~\cite{harris2020}, but numerous factors influence their complex coupling, and the underlying mechanisms are still elusive. 
Also, there is an intrinsic difficulty of isolating these factors in a clinical setting and measuring their contribution~\cite{vercellin2016}. 
The present work focuses on the interaction between the main ocular vessels' hemodynamics, intraocular pressure and the retrolaminar tissue pressure, which is directly related to the cerebrospinal fluid pressure. 
Among several interesting contributions in this area, we mention those closely related to our work. 
The first mathematical model that simultaneously accounts for blood flow in the central retinal vessels, blood flow in the retinal microvasculature, retinal blood flow autoregulation, biomechanical action of intraocular pressure on the retinal vasculature, and time-dependent arterial blood
pressure was introduced in~\cite{guidoboni2014intraocular}. 
A theoretical model to study the effects of intraocular pressure elevation on the central retinal artery hemodynamics was proposed in ~\cite{guidoboni2014effect} and extended to account for the central retinal venous hemodynamics and the retinal microcirculation in~\cite{carichino2014}.
However, none of them explicitly accounted for uncertainties and variabilities in the model parameters. Only a few modelling works include a stochastic analysis framework and focused on the production and drainage of aqueous humour flow~\cite{szopos2016} and its coupling with ocular hemodynamics in a simplified manner~\cite{sacco2015}.\\

With these premises, we present in this contribution an uncertainty quantification and a global sensitivity analysis for the main parameters involved in the mathematical and computational framework called the Ocular Mathematical Virtual Simulator (OMVS) that we have developed~\cite{sala2018ocular,sala2019mathematical}. 
The clinical relevance of results provided by the OMVS is described in~\cite{sala2017patient} or, more extensively, in~\cite[Chapter 13]{sala2019mathematical}, and it has been confirmed by an independent population-based study including nearly $10000$ individuals~\cite{tham2018inter}. 
The reduced version of the OMVS model employed in the present study originates from~\cite{guidoboni2014intraocular} for the retinal circulation.
It has been extended to include a reduced model for blood flow perfusion in the lamina cribrosa. Preliminary findings from a simplified uncertain quantification study were published as a peer-reviewed conference abstract in~\cite{sala2019cmbe}. As a significant step forward, we present hereafter a detailed global sensitivity analysis, using the Sobol' sensitivity indices.\\

The paper is organized as follows. The mathematical and computational model is described
in Section~\ref{subsec:math_comp}, the uncertainty quantification and sensitivity analysis approach 
is presented in Section~\ref{subsec:uq_sa} and the input data in Section~\ref{subsec:input_data}. The results of our study are depicted in Section~\ref{sec:results} and discussed in Section~\ref{sec:discussion}. Finally, conclusions and future perspectives are outlined in Section~\ref{sec:concl}.

\section{Methodology} \label{sec:methods}

In the next two sections we describe the deterministic mathematical and computational foundations of our study, as well as the uncertainty quantification (UQ) and sensitivity analysis (SA) methods we incorporated in the OMVS framework to account for the stochastic features of the system.

\subsection{Mathematical and computational model} \label{subsec:math_comp}
 The OMVS is a complex modelling framework that couples hemodynamics, biomechanics, and fluid dynamics in the eye to visualize and estimate in a non-invasive way ocular biofluids and tissues characteristics that are difficult or not accessible with standard investigation methods. 
 The contributions developed within this framework can be subsequently utilized to isolate single risk factors and quantify their influence on the multi-factorial disease process. \\

To achieve this goal, the full OMVS is designed with a multiscale architecture, that aims at preserving the natural systemic features of blood circulation, while providing detailed views on sites of particular interest from the clinical viewpoint, such as the lamina cribrosa. 
The lamina cribrosa is a sponge tissue in the back of the eye that has a crucial role from the hemodynamical and neurological viewpoints. 
Specifically, this membrane is thought to help maintain the balance between the pressure inside the eye (intraocular pressure, hereafter denoted IOP) and behind the eye (the retrolaminar tissue pressure directly influenced by the intracranial pressure, denoted RLTp), which may influence the ocular blood flow. 
In addition, the lamina cribrosa acts as a scaffold for the retinal ganglion cell axons and the central retinal vessels and feeds RGC axons through its vascular network. 
The IOP is easily measurable with a Goldmann applanation tonometer. 
This instrument is based on the Imbert–Fick principle, which affirms that the pressure inside a dry thin-walled sphere corresponds to the force required to flatten the sphere surface divided by the flattening area.

\begin{figure}[h!]
	\centering
	\begin{subfigure}{0.495\linewidth}
		\centering
		\includegraphics[width=\linewidth]{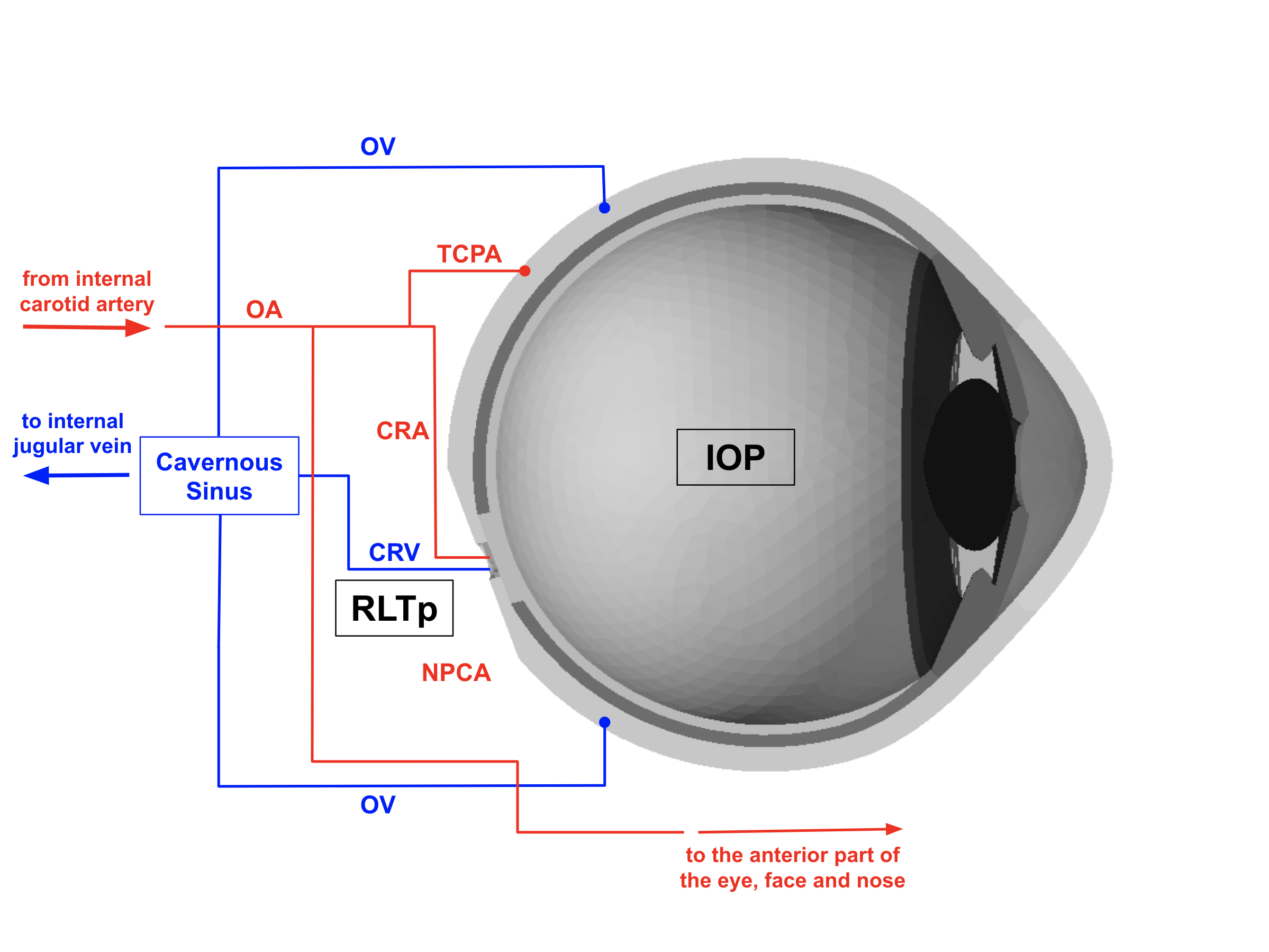}
		\caption{System I: 0D hemodynamics of retinal vasculature.}
		\label{fig:omvs_scheme:systemI}
	\end{subfigure}
	\begin{subfigure}{0.495\linewidth}
		\centering
		\includegraphics[width=\linewidth]{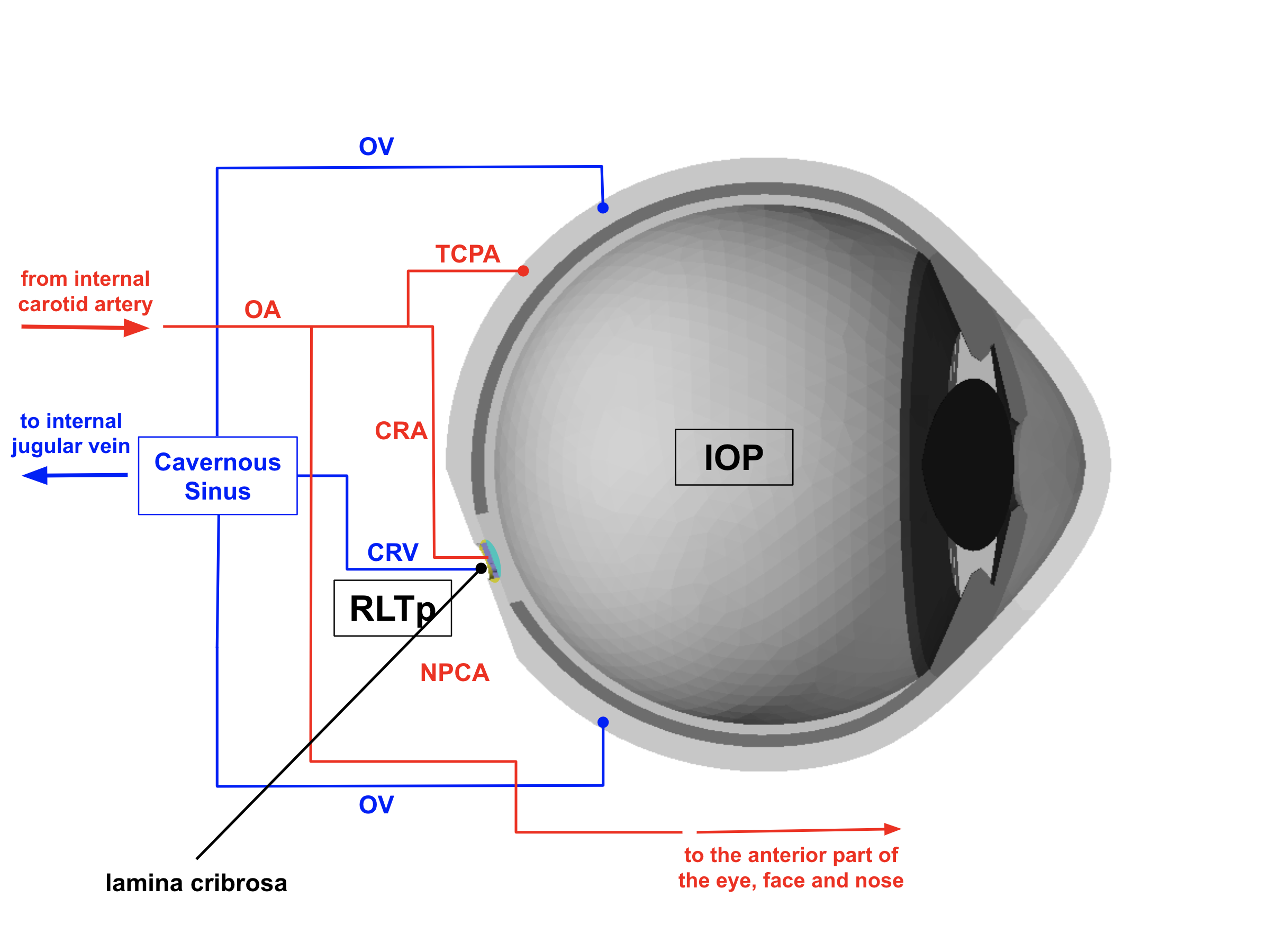}
		\caption{System II{: 0D hemodynamics of retinal vasculature + 3D hemodynamics of lamina cribrosa.}}
		\label{fig:omvs_scheme:systemII}
	\end{subfigure}
	\begin{subfigure}{0.495\linewidth}
		\centering
		\includegraphics[width=\linewidth]{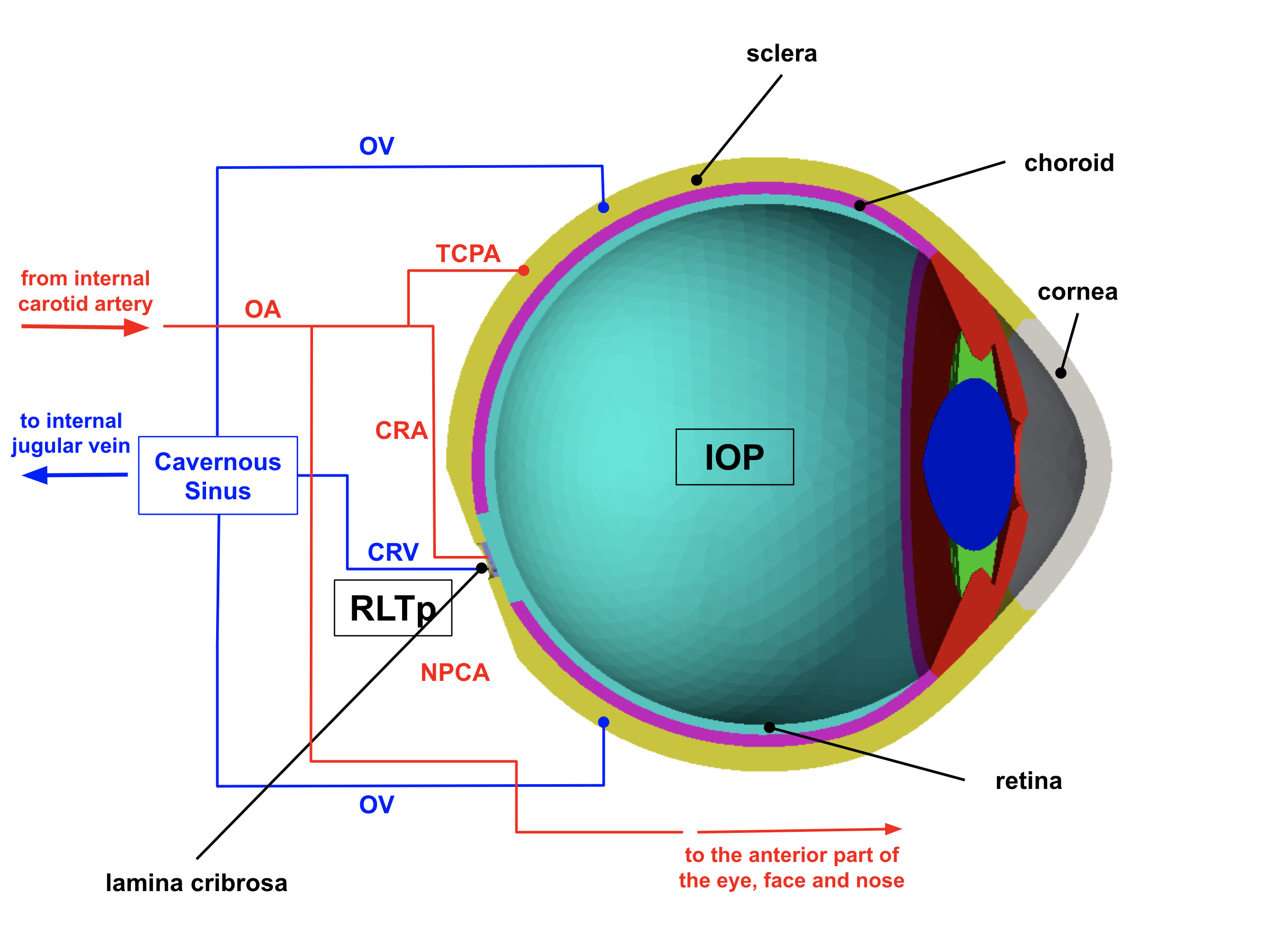}
		\caption{System III{: 0D hemodynamics of retinal vasculature + 3D hemodynamics of lamina cribrosa + 3D biomechanics of lamina cribrosa, retina, choroid, sclera, and cornea.}}
		\label{fig:omvs_scheme:systemIII}
	\end{subfigure}
	\caption{Schematic of the OMVS multiscale structure. 
		IOP = intraocular pressure, RLTp = retrolaminar tissue pressure, CRA = Central Retinal Artery, CRV = Central Retinal Vein, OA = ophthalmic artery, OV = ophthalmic vein, TCPA = temporal posterior ciliary artery, NPCA = nasal posterior ciliary artery.}
	\label{fig:omvs_scheme}
\end{figure}

In practice,  we developed three model formulations of increasing complexity: we started from a $0$-dimensional reduced-order description of the system, progressively adding the coupling with a porous media model for the lamina cribrosa and finally incorporating the effects of the deformation of the ocular tissues. 
More precisely, the OMVS combines (see Fig.~\ref{fig:omvs_scheme}):

\begin{enumerate}
	\item System I (Fig.~\ref{fig:omvs_scheme:systemI}): a circuit-based ($0$D) model for blood flow in the retinal vasculature, central retinal artery (CRA), and central retinal vein (CRV);
	\item System II (Fig.~\ref{fig:omvs_scheme:systemII}): a three-dimensional ($3$D) porous media model for the perfusion of the lamina cribrosa;
	\item System III (Fig.~\ref{fig:omvs_scheme:systemIII}): a $3$D isotropic elastic model for the biomechanics of the lamina cribrosa, retina, choroid, sclera, and cornea.
\end{enumerate}

In the present contribution, UQ and SA analyses require intensive evaluations of the physical-based model. 
Therefore, we have employed a reduced version of the OMVS, accounting for the hemodynamical description provided by System I (Fig.~\ref{fig:omvs_scheme:systemI}), coupled to a $0$D adaptation of System II (Fig.~\ref{fig:omvs_scheme:systemII}). 
For the proposed study, System III (Fig.~\ref{fig:omvs_scheme:systemIII}) has not been considered. 
In this manner, the model (i) provides a multiscale hemodynamics overview of the overall system, while maintaining a relatively accessible mathematical complexity and low computational costs; and (ii) combines information on ocular sites for which quantitative data are available - \textit{e.g.} blood flow in the central retinal artery - and crucial ocular areas that are not accessible with clinical images - \textit{e.g.} lamina cribrosa perfusion.

\begin{figure}[h!]
	\centering
	\includegraphics[width=0.9\linewidth]{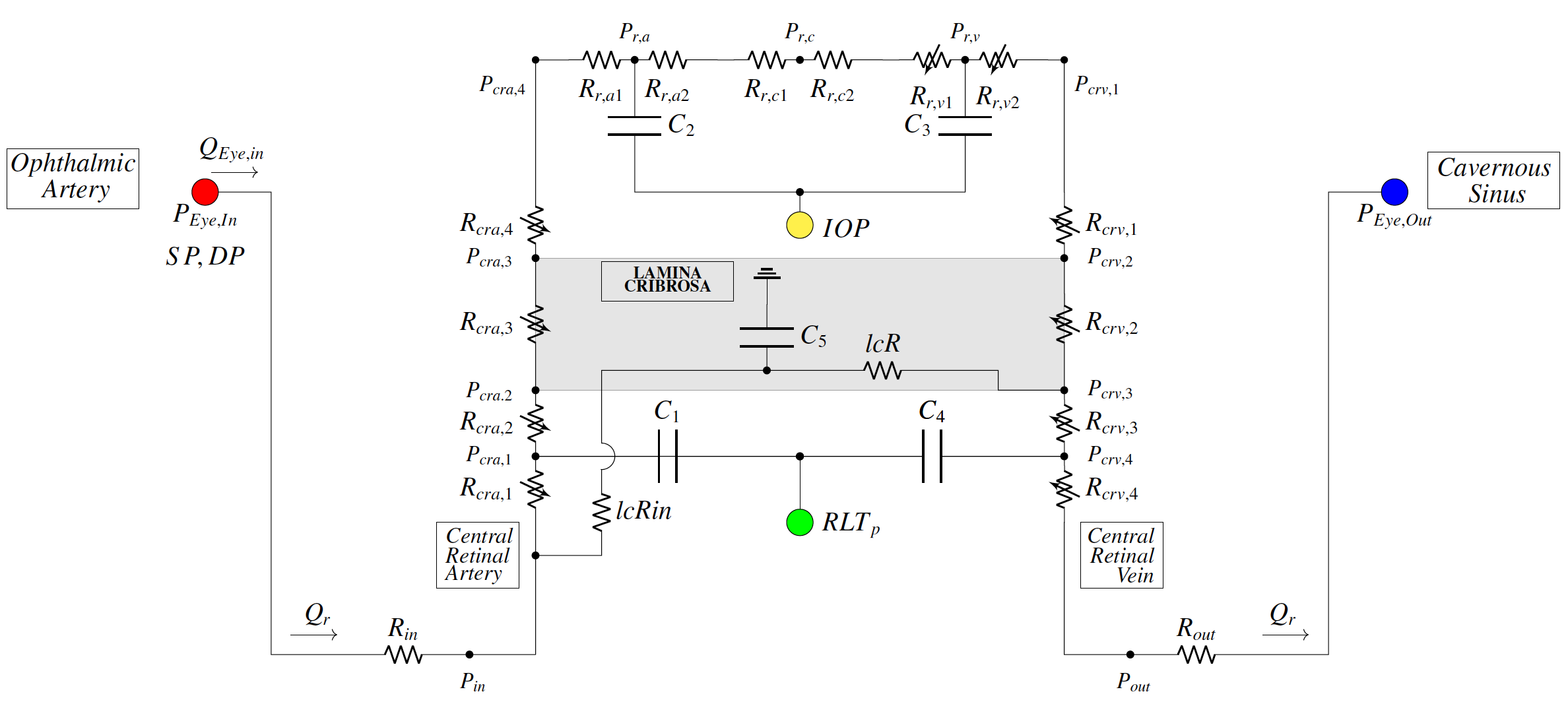}
	\caption{Schematic of the $0$D circuit-based model {(System I of the OMVS) employed in this contribution} for UQ and SA. 
		The coloured nodes are the input of the model that can be inferred from clinical measurements.}
	\label{fig:circuit_scheme}	
\end{figure}

The $0$D reduced version of the OMVS, see Fig.~\ref{fig:circuit_scheme}, exploits the electric analogy to fluid flow in complex vascular network~\cite{formaggia2010cardiovascular}. 
In this context, electric potentials correspond to fluid pressure, electric charges correspond to fluid volumes, and electric currents correspond to volumetric flow rates; the resistors and capacitors represent hydraulic resistance and wall compliance, respectively. 
Writing the constitutive equations characterizing the circuit elements and the Kirchhoff laws of currents and voltages leads to a system of ordinary differential equations whose solution provides the time-dependent profiles of pressures at the circuit nodes and flow rates through the circuit branches.
We emphasize that IOP plays a crucial role in the description of the vein collapsibility, which is modelled in the OMVS by Starling resistors~\cite{starling1896absorption}.
Namely, when the external pressure is higher than the internal blood pressure, veins collapse, therefore dropping down the blood flow.

The network is constructed as an extension of a previous model for the retinal circulation, proposed and validated in~\cite{guidoboni2014intraocular} and ~\cite{cassani2016blood}. 
The vasculature is divided into six main compartments: central retinal artery (\textit{cra}), arterioles (\textit{r,a}), capillaries (\textit{r,c}), venules (\textit{r,v}), central retinal vein (\textit{crv}), and the lamina cribrosa (\textit{lc}). 
Each compartment includes resistances (R) and capacitances (C). 
The intraocular segments ($R_{cra,3},R_{cra,4}, R_{r,v1},R_{r,v2},R_{crv,1},R_{crv,2}$) are exposed to the IOP and the retrobulbar segments are exposed to the RLTp ($R_{cra,1},R_{cra,2}, R_{crv,3},R_{crv,4}$). 
The explicit conservation and constitutive laws, as well as parameters involved in the description of the retinal circuit, follow directly from work by ~\cite{guidoboni2014intraocular} and ~\cite{cassani2016blood}.
To the initial model, we have added a simplified description of the hemodynamics in the lamina cribrosa, involving two resistors $lcRin$ and $lcR$, and one capacitor $C_5$, with the following values: $lcRin=78181.9$ \unit{mmHg s cm}$^{-3}$, $lcR=23988.25$ \unit{mmHg s cm}$^{-3}$, and $C_5 = 0.000000753$ \unit{cm}$^3$ \unit{mmHg}$^{-1}$.
Also, in the original circuit, the external pressure on the resistances $R_{cra,3}$ and $R_{crv,2}$ is the effective stress exerted by the lamina on these vessels, which has been computed \emph{via} a simplified  fluid-structure interaction model~\cite{guidoboni2014effect} describing the CRA/CRV interaction with the lamina cribrosa. 
In the current version of the model, we do not account for this contribution in a similar manner but rather adopt a simplified approach, in which the external pressure corresponds to IOP. Further extensions could incorporate this dependence, but our choice was dictated by the possibility of computing stresses directly from the System III component of the OMVS. \\
 
The reduced model thus obtained can predict the hemodynamics within the lamina cribrosa, the retinal vasculature, and the central retinal vessels based on the key inputs described in Tab.~\ref{tab:input_params} and displayed as coloured dots in Fig.~\ref{fig:circuit_scheme}.
For the proposed study we fixed the value of the pressure at cavernous sinus (blue node), while we choose as input variables for UQ and SA the systolic and diastolic pressure at the ophthalmic artery (SP and DP, red node), the intraocular pressure (IOP, yellow node) and the retrolaminar tissue pressure (RLTp, green node).

Regarding the outputs, and in light of the clinical application in view, we will focus on the quantities of interest in the UQ and SA analysis listed in Tab.~\ref{tab:output_params}. Note that these quantities of interest will reflect the behaviour of the system at different time instants through the cardiac cycle.

 \begin{table}[h!]
 	\centering
	\begin{tabular}{r c l}
 		\textbf{Key name} & \textbf{Unit} & \textbf{Brief description}  \\ 
 		\hline \noalign{\smallskip}
 		DP & \unit{mmHg} & Diastolic blood Pressure  \\
 		SP & \unit{mmHg} & Systolic blood Pressure \\
 		IOP & \unit{mmHg} & IntraOcular Pressure \\
 		RLTp & \unit{mmHg} & RetroLaminar Tissue pressure \\
 		\hline \noalign{\smallskip}
 	\end{tabular}
 	\caption{Input parameters for the reduced OMVS.}
 	\label{tab:input_params}
	 \end{table}	
 	\begin{table}[h!]
 		\centering
 		\begin{tabular}{r c l}
 				\textbf{Key name} & \textbf{Unit} & \textbf{Brief description}  \\ 
 				\hline \noalign{\smallskip}
 				CRA\_ps & $\unit{\mu l/min}$ & peak systolic CRA blood flow  \\
 				CRA\_es & $\unit{\mu l/min}$ & end systolic CRA blood flow  \\
 				CRA\_ed & $\unit{\mu l/min}$ & end diastolic CRA blood flow  \\
 				CRV\_ps & $\unit{\mu l/min}$ & peak systolic CRV blood flow  \\
 				CRV\_es & $\unit{\mu l/min}$ & end systolic CRV blood flow  \\
 				CRV\_ed & $\unit{\mu l/min}$ & end diastolic CRV blood flow  \\
 				LC\_ps & $\unit{\mu l/min}$ & peak systolic lamina cribrosa blood flow  \\
 				LC\_es & $\unit{\mu l/min}$ & end systolic lamina cribrosa blood flow  \\
 				LC\_ed & $\unit{\mu l/min}$ & end diastolic lamina cribrosa blood flow  \\
 				\hline \noalign{\smallskip}
 		\end{tabular}
 		\caption{Quantities of interest in the UQ and SA analysis.}
 		\label{tab:output_params}
 \end{table}

\begin{figure}[h!]
	\centering
	\includegraphics[width=0.7\linewidth]{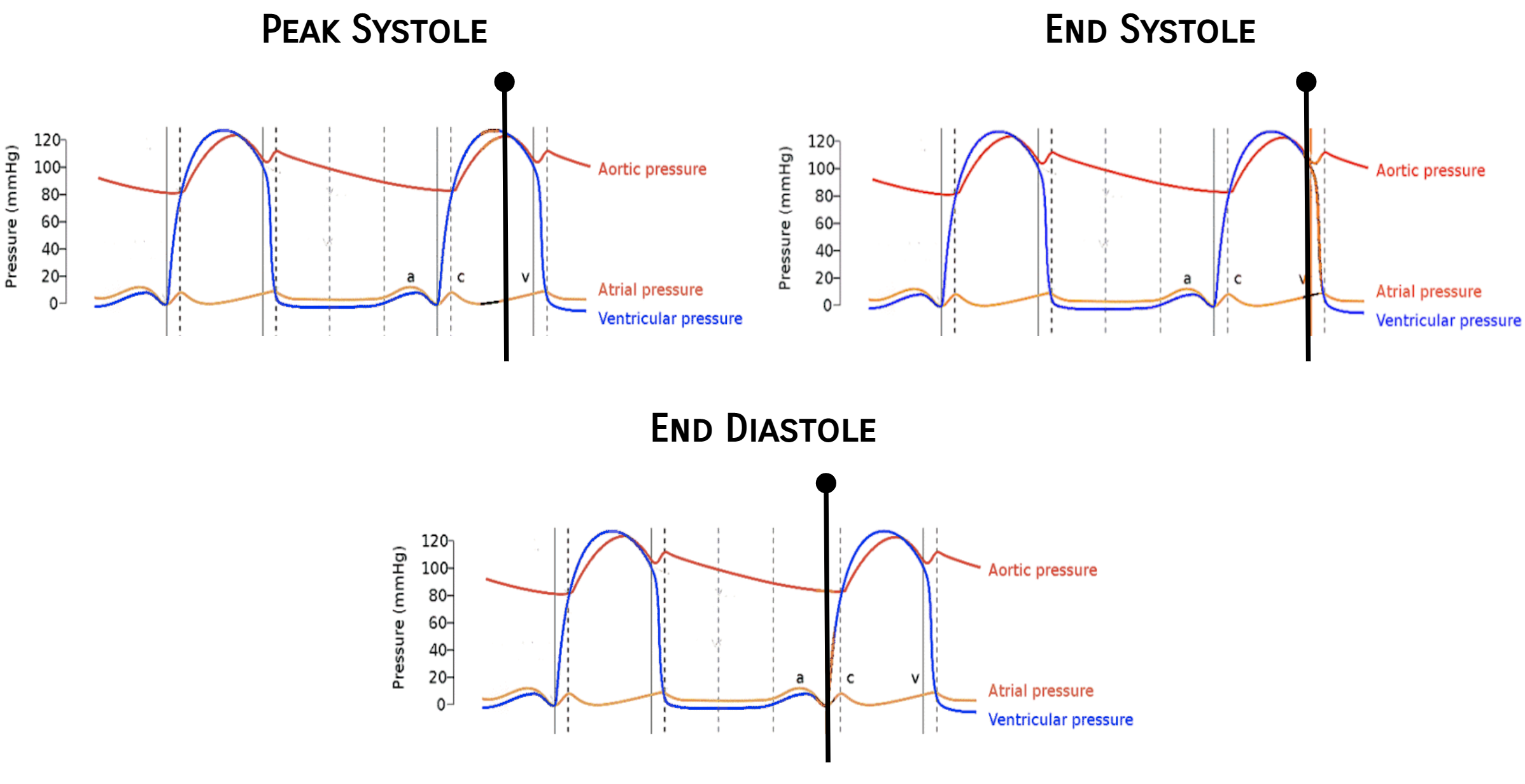}
	\caption{{Specific instant during the last cardiac cycle (peak systole, end systole and end diastole) used to compute the clinical outputs of interest. Image edited from \url{https://en.wikipedia.org/wiki/Cardiac_cycle}.}}
	\label{fig:cardiac_cycle}
\end{figure}

The mathematical model previously described has been implemented in OpenModelica~\cite{fritzson2006openmodelica}, an open-source Modelica-based modelling and simulation environment intended for industrial and academic studies of complex dynamic systems. 
Model results have been obtained using DASSL~\cite{petzold1982description} with a tolerance of $10^{-6}$, a time step of $10^{-3}$ and total simulation time of $8 \, s$.
DASSL (Differential/Algebraic System Solver) is an implicit, high order, multi-step solver with a step-size control based on backward differentiation formula (BDF).
These features allow it to be stable and fit to be used for a wide range of models. Its first development can be in found in \cite{petzold1982description}. \\
The system reaches a periodic state after the first cardiac cycle. 
However, we consider our output the last simulated cardiac cycle. 
Then, we retrieve the CRA, CRV, and the lamina cribrosa blood flow at three specific instant during the last cardiac cycle, namely the peak systolic time, the end of the systole and the end of the diastole {(see Fig. \ref{fig:cardiac_cycle})}. 
{The choice of these three particular time instant in the cardiac cycle is driven by their interest from a clinical perspective. Moreover, in view of the model validation, the measurements of the blood flow is very often taken at peak systole and end diastole \cite{harris1996acute}.}

\subsection{Uncertainty quantification and sensitivity analysis approach}
\label{subsec:uq_sa}

\paragraph{Uncertainty propagation.} The construction of a reliable model for ocular biofluid dynamics involves several steps with inherent uncertainties, among which: parameter inference from uncertain experimental data, model personalization to the same subject data at different time instants or to different individuals, {\it etc.} 
Therefore, a major challenge is to assess how these sources of uncertainty impact the clinically relevant outputs of the simulation and ultimately affect the confidence in the model predictions.\\

In the present contribution, we adopted the following approach: for the set of key inputs of the reduced version of the OMVS model, the uncertainty is represented by a probability density function (pdf), which quantifies the probability of a given parameter to reproduce a specific observation. We next develop a forward uncertainty quantification, also known as uncertainty propagation, to investigate how these input uncertainties are propagated to the outputs {\it via} the computational model. Thus, the combined information between previous modelling knowledge and assumptions on the prior pdfs allows us to obtain the posterior pdf for the quantities of interest selected from the clinical perspective.\\

\begin{remark}
A major challenge that needs to be addressed is the choice and representation of the prior pdf in light of the quantitative and qualitative information available from generally noisy data, see for instance~\cite{SundnesCMBE2019} and further discussion in Sec.~\ref{subsec:input_data}.
\end{remark}

\paragraph{Sensitivity analysis.} This part aims to determine {synthetic measurements of}  which key inputs are the most influential 
on the quantities of interest selected among the outputs of the computational model {without making assumption on the model and taking into account the continuous nature of the input parameters}. To this end, we adopt 
the stochastic framework of global sensitivity analysis, which considers the input parameters $\{X_j\}_{j\in \{1, \ldots d\}}$ to be random independent variables with uncertainty modelled by a probability distribution, and employed to compute the random output $Y$. We have not considered the sensitivity analysis of the dynamic process $Y(t)$. 
The methods described hereafter can be adapted to time $t$, and we could study the sensitivity indices with respect to time. 
However, it is more useful {--- easier to interpret ---} to perform sensitivity analysis at important characteristics  of the  response time series, in our case, peak systole, end systole, and end diastole, as discussed in Sec.~\ref{subsec:math_comp}.

To quantify the influence of the variations of $X_j$ on the variations of $Y$, we compute the so-called Sobol' sensitivity indices originally proposed in the seminal paper~\cite{sobol1993sensitivity}, see also~\cite{prieur2017variance}. More precisely, we define the first-order indices as
\begin{equation}
 \label{eq:sobol:1st}
  S_j = \frac{\var[\mathbb{E}(Y|X_j)]}{\var[Y]},
\end{equation}
where 
\begin{inparaenum}[(\it i)]
\item $\var$ denotes the variance and $\mathbb{E}$ the expected value;
\item $\var[Y]$ corresponds to the  variability of $Y$ with respect to the overall uncertainty including non-linear effects;
\item $\var[\mathbb{E}(Y|X_j)]$, the variance of the conditional expectation $\mathbb{E}(Y|X_j)$, corresponds to the  main or first order  effect of $X_j$; it means that if $Y$ is sensitive to $X_j$, $\mathbb{E}(Y|X_j)$ is likely to vary a lot and hence $\var[\mathbb{E}(Y|X_j)]$ as well.
\end{inparaenum}

Another useful index is the total Sobol' index, defined as 
\begin{equation}
 \label{eq:sobol:total}
  S_j^{tot} = 1 - \frac{\var[\mathbb{E}(Y|X_{(-j)})]}{\var[Y]} = 1 - S_{-j}, 
\end{equation}
where $X_{(-j)} = (X_1,\ldots,X_{j-1},X_{j+1},\ldots,X_d)$ and $S_{-j}$ denotes the sum of the indices where $X_j$ is not involved.
Additionally, the effect due to specific interactions between the $j^{th}$ and the $k^{th}$ factors ($k\neq j$) can be measured by second-order Sobol' indices, and so on for high-order interactions, see for more details~\cite{prieur2017variance}.

Several approaches have been proposed to numerically compute these sensitivity indices, as reviewed for instance in~\cite{prieur2017variance}. 
In the present work, we adopted the following two strategies: (i) a Monte Carlo-type approach and an estimator proposed in~\cite{saltelli2002} on the basis of a combinatoric argument, and (ii) a Fourier amplitude sensitivity test (FAST)~\cite{saltelli1999}, which is a spectral method based upon the Fourier decomposition of the model response. 
The computational cost of the first approach for first order and total order Sobol' indices is of $(d+2)n$ model evaluations~\cite{saltelli2002}, where $d$ is the input space dimension and $n$ is the sample size; as for the FAST method the computation cost for first and total order indices is of $d \; n$ model evaluations~\cite{saltelli1999}.

This variance-based approach implies intensive sampling, but it allows us to explore the input factors' full uncertainty ranges. We provide the two strategies as a way to ensure  the reliability of our estimates. Indeed, given a sampling size,  they may  vary  when repeating the estimations, and, in the Monte-Carlo approach, they can even be negative, although we used an implementation that mitigates this issue.

All the results on the UQ and SA analysis presented hereafter are carried out exploiting the Python statistical library OpenTURNS~\cite{baudin2017openturns}.\\

\subsection{Input data} \label{subsec:input_data}
Both mathematical methods described before need as prior knowledge the statistical distribution of the input.
We detail and critically discuss in the sequel the choices we propose, based on assumptions deduced from the experimental and clinical literature.
Recall that the input of our model are systolic blood pressure (SP), diastolic blood pressure (DP), IOP and RLTp.

\begin{figure}[h!]
	\centering
	\begin{subfigure}{0.32\linewidth}
		\centering		
		\includegraphics[width=\linewidth]{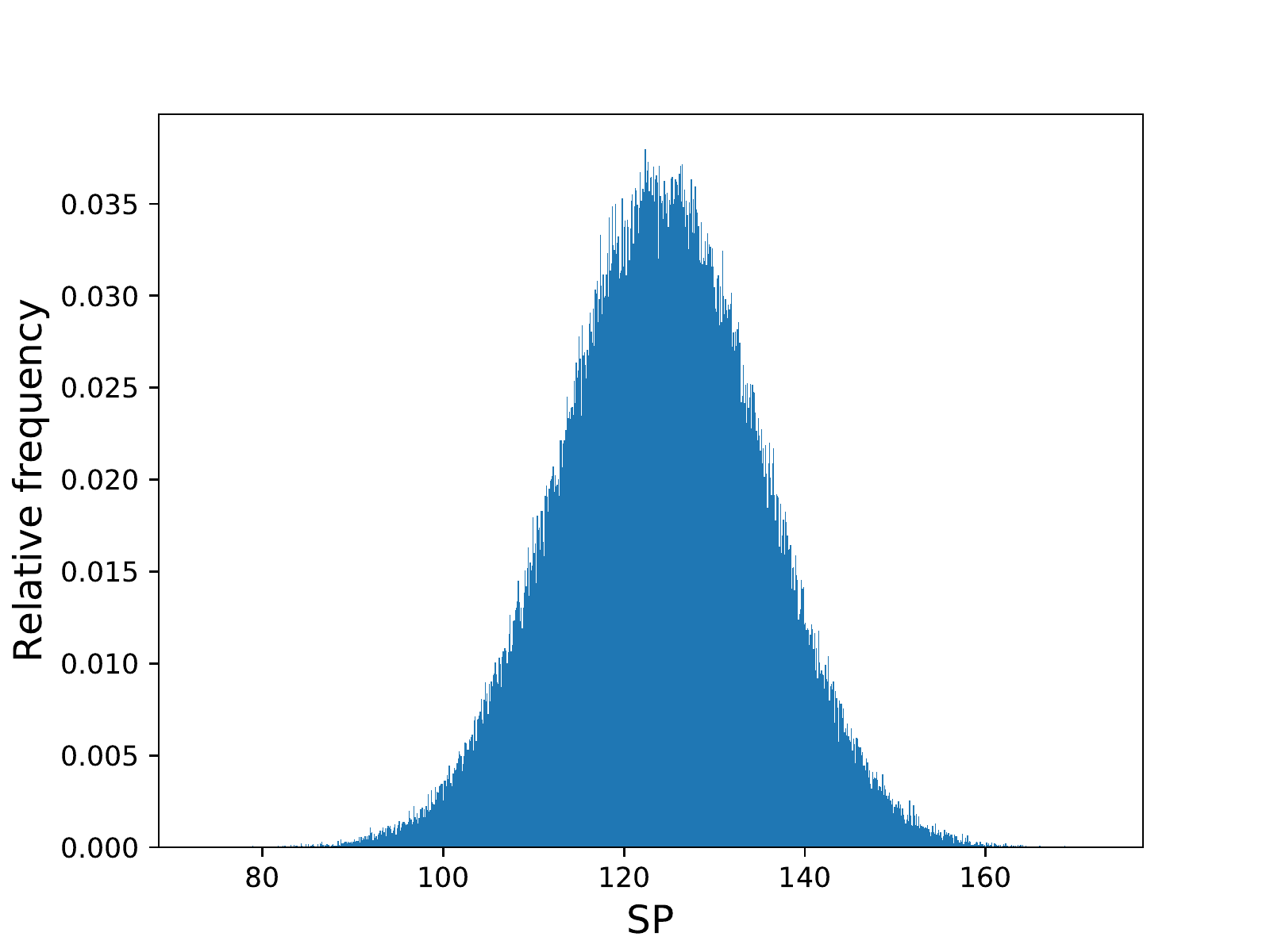}
		\caption{Systemic systolic blood pressure}
		\label{fig:input:pdf_sp}
	\end{subfigure}
	\begin{subfigure}{0.32\linewidth}
		\centering
		\includegraphics[width=\linewidth]{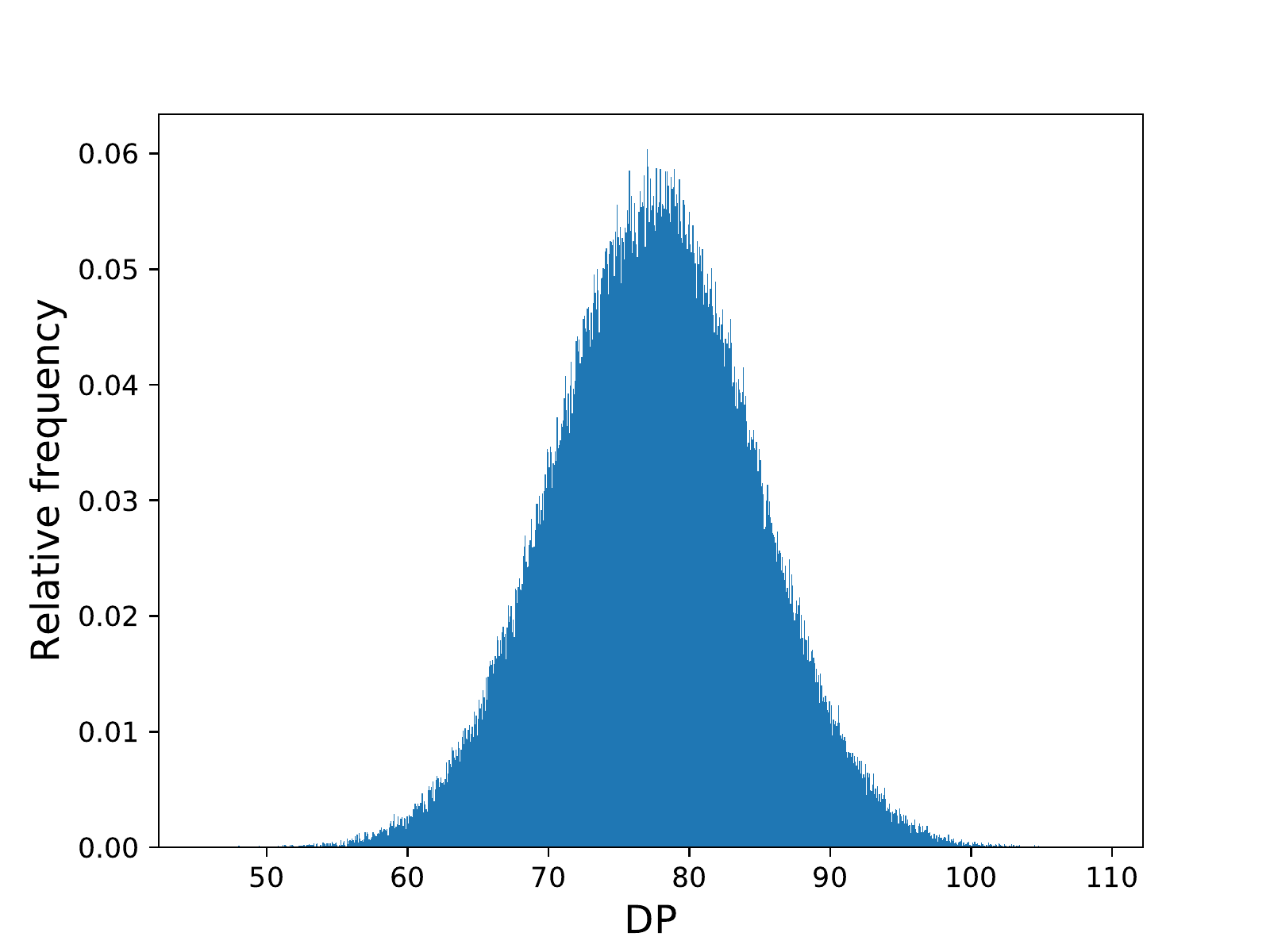}
		\caption{Systemic diastolic blood pressure}
		\label{fig:input:pdf_dp}
	\end{subfigure}
	\begin{subfigure}{0.32\linewidth}
		\centering
		\includegraphics[width=\linewidth]{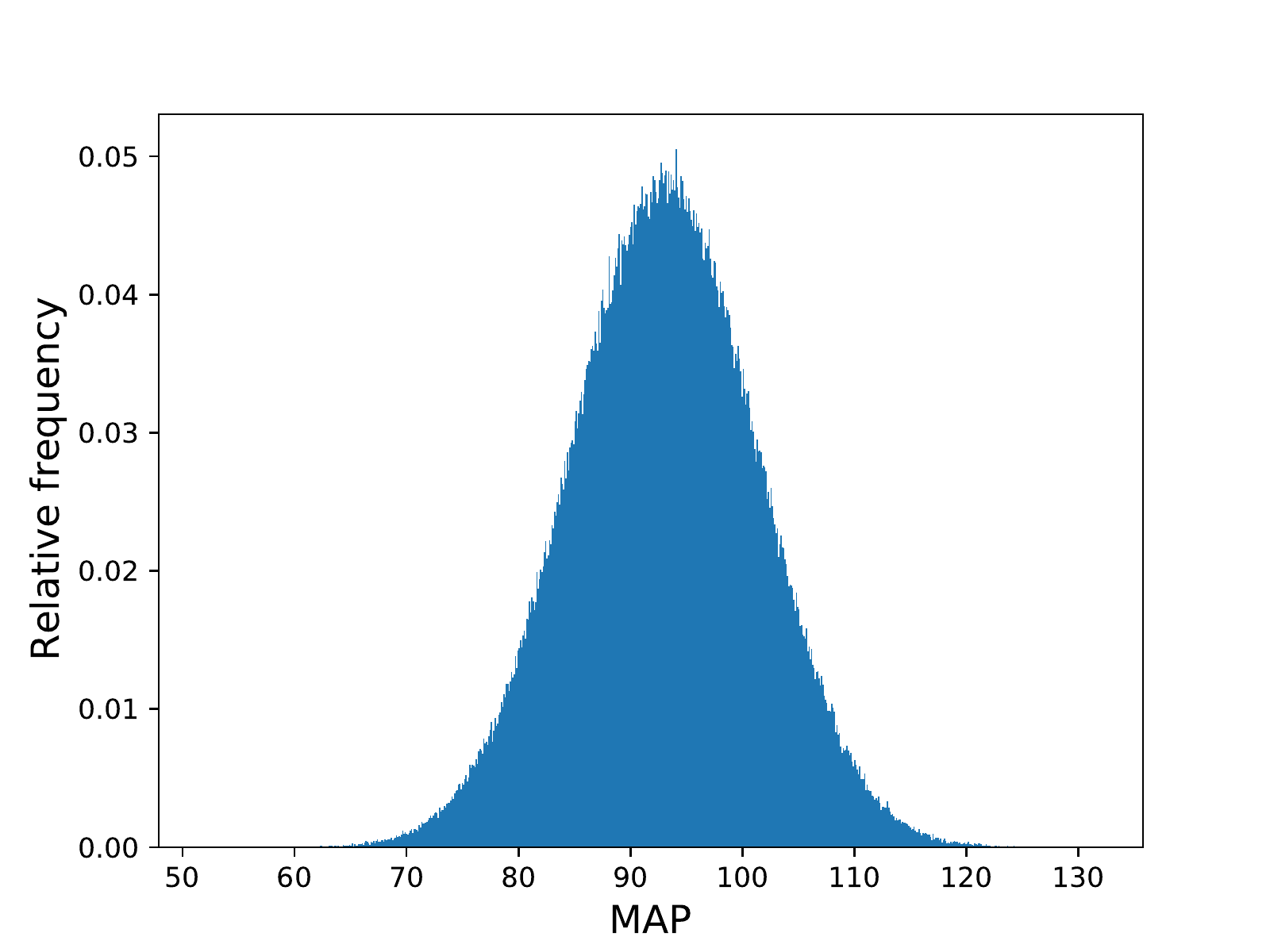}
		\caption{Systemic mean arterial pressure}
		\label{fig:input:pdf_map}
	\end{subfigure}
	\caption{{Normal} probability density functions of the input blood pressure variables}
\end{figure}

\paragraph{Blood pressures.} 
For SP and DP data we refer to the paper of Sesso and co-authors~\cite{sesso2000systolic} where these two quantities showed a normal distribution.
In particular SP has a mean of $124.1 \, \unit{mmHg}$ and a standard deviation of $11.1 \, \unit{mmHg}$  (Fig. \ref{fig:input:pdf_sp}), whereas DP has a mean of $77.5 \, \unit{mmHg}$ and a standard deviation of $7.1 \, \unit{mmHg}$ (Fig.~\ref{fig:input:pdf_dp}).

Using SP and DP as input parameters raises an issue for the sensitivity analysis, notably the Sobol' index study, namely that the inputs have to be assumed independent, see Sec.~\ref{subsec:uq_sa}.
For IOP, RLTp, we can make this assumption reasonably, however, this is not valid for SP and DP \cite{gavish2008linear}.
To overcome this problem we used Mean Arterial Pressure (MAP) as input variable in the Sobol study that has a normal distribution with mean $93 \, \unit{mmHg}$ and standard deviation of $7.6 \, \unit{mmHg}$.
In this case we also reconstructed the MAP distribution starting from the normal distribution of SP and DP and the following relationship:
\begin{equation} \label{eq:map}
	MAP = \frac{1}{3} SP + \frac{2}{3} DP
\end{equation}
to check its normal distribution assumption (Fig.~\ref{fig:input:pdf_map}). \\
From Eq.~\eqref{eq:map} and from the correlation assumption between SP and DP made in~\cite{gavish2008linear}, we can then reconstruct SP and DP starting from the MAP.
Finally, the MAP is independent of the other two inputs (IOP, RLTp) and can be used as a sensitivity analysis parameter.

\paragraph{Intraocular pressure.}
There is a considerable discussion about normal and lognormal probability density functions in literature, and which of the two can better represent biological phenomena~\cite{limpert2001log}.

To explain our modelling choices we consider different IOP distributions using the data recovered from the clinical work of Suh and collaborators~\cite{suh2012distribution} (mean $\mu_{normal}=14.7 \,\unit{mmHg}$, standard deviation $\sigma_{normal}=2.8\, \unit{mmHg}$).
Starting from these values, we have computed the mean and the variance both for Gaussian distribution and for a lognormal behaviour using the following formulas~\cite{johnson1994lognormal} to be consistent with the data:
\begin{align*}
\mu_{lognormal} &= \log \left( \frac{\mu_{normal} }{\sqrt{1+\frac{\sigma^2_{normal}}{\mu^2_{normal}}}}\right) &
\sigma^2_{lognormal} &= \log \left(  1+ \frac{\sigma^2_{normal}}{\mu^2_{normal}}\right)
\end{align*}

We performed a comparative analysis using three different IOP distributions based on the same clinical data, which have been defined above. In particular, we compare a \textit{normal}, a \textit{truncated normal} and a \textit{lognormal} distribution.

\begin{figure}[h!]
	\centering
	\begin{minipage}{0.49\linewidth}
	\centering
	\includegraphics[width=\linewidth]{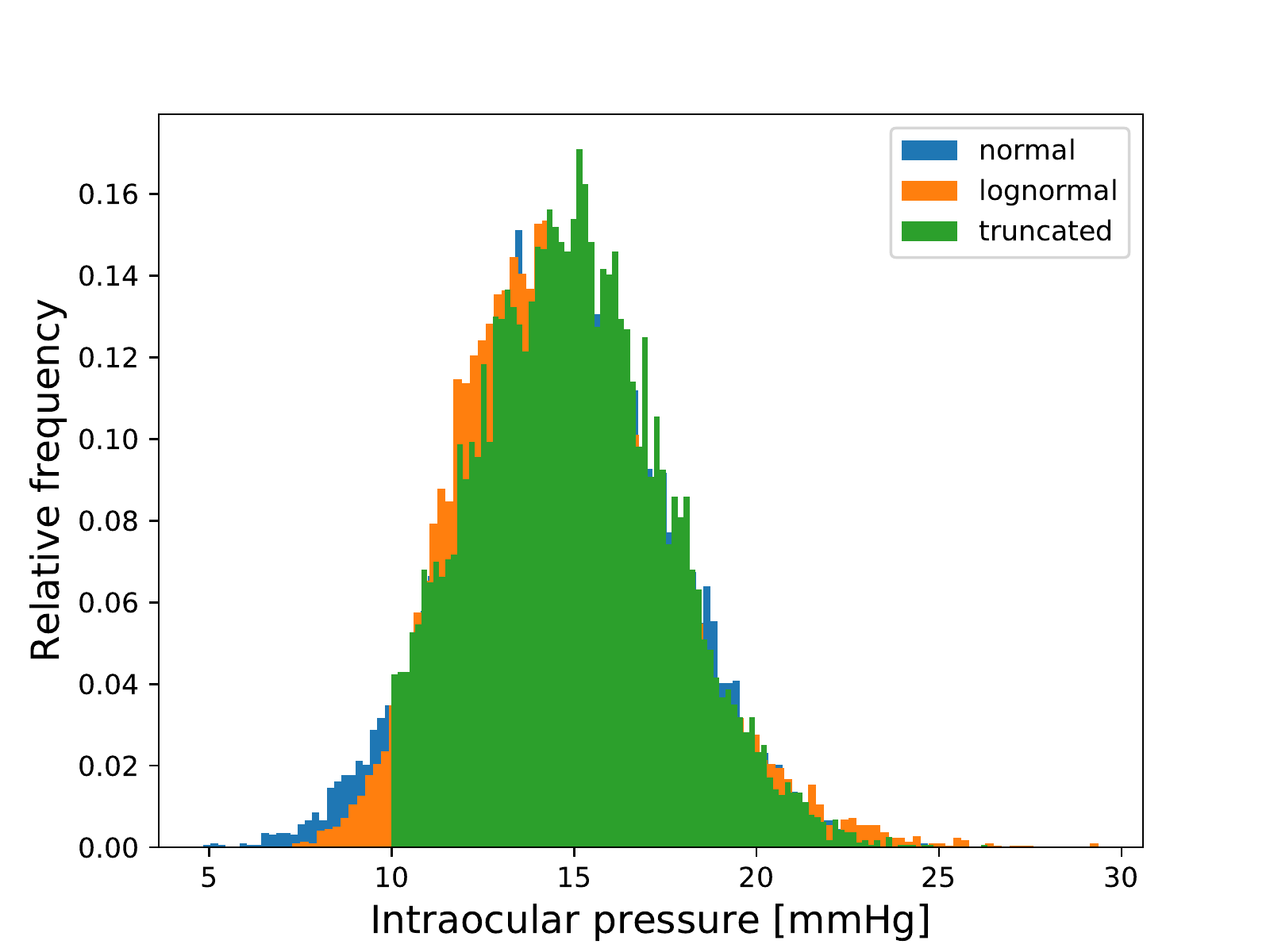}
	\caption{Comparison among normal (blue), lognormal (orange) and truncated normal (green) IOP distributions using the same input data \cite{suh2012distribution}.}
	\label{fig:input:comparison-IOP}
	\end{minipage}
	\begin{minipage}{0.49\linewidth}
	\centering
	\includegraphics[width=\linewidth]{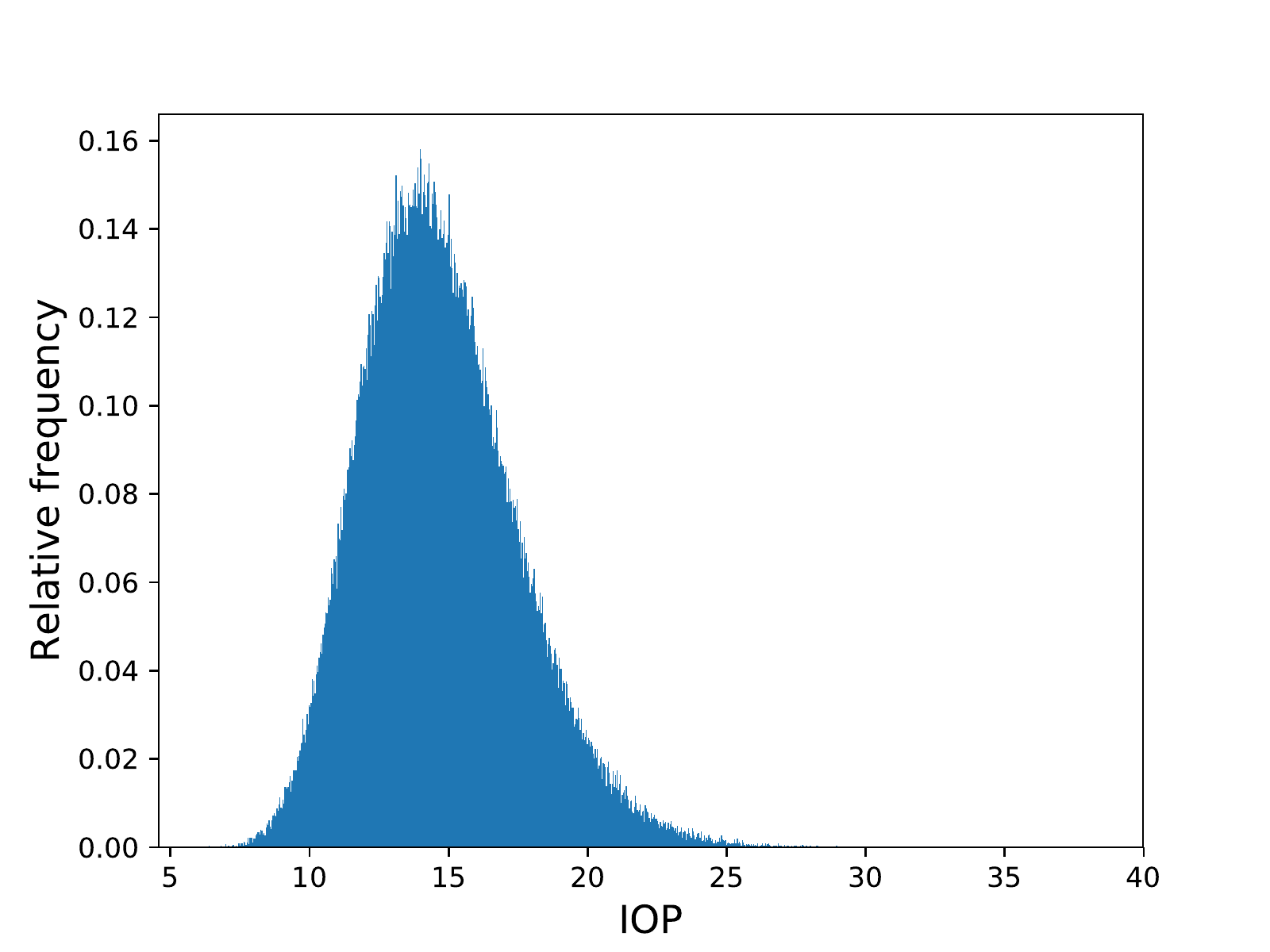}
	\caption{{Lognormal} probability density function of the input intraocular pressure (IOP).}
	\label{fig:input:pdf_iop}
	\end{minipage}
\end{figure}

Fig.~\ref{fig:input:comparison-IOP} highlights that the normal probability density (blue) is going beyond some physiological constraints for healthy patients such as $IOP > 5 \unit{mmHg}$.
The truncated normal probability density function (green) does not show this issue; however, this IOP distribution presents an abnormal cut in the left tail, making it ineffective for a sensitivity analysis study. For our simulations, with the data provided by Suh et al.~\cite{suh2012distribution}, the lognormal distribution seems to be the more natural one.
Our choice is dictated by the fact that we want to avoid miscalculation due to unphysiological input parameters that may lead to unrealistic discontinuities in the simulation results.
This lognormal assumption has also been accepted in other works~\cite{wang2011intraocular}.
For this reason we utilized the log normal distribution based on the population based study operated by Suh et al. ~\cite{suh2012distribution} (Fig.~\ref{fig:input:pdf_iop}).

\begin{figure}[h!]
	\centering
	\includegraphics[width=0.5\linewidth]{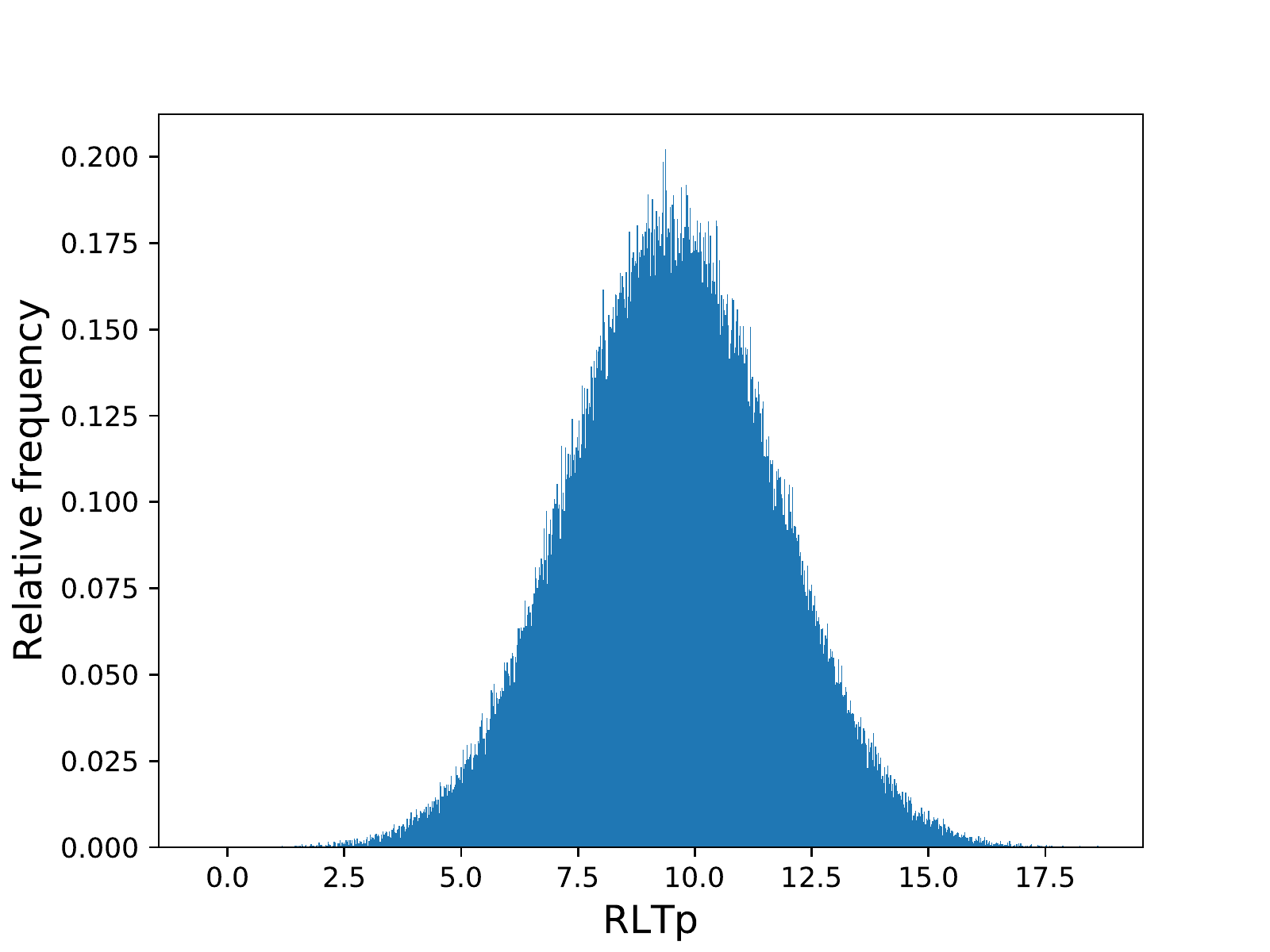}
	\caption{{Normal} probability density function of the input retrolaminar tissue pressure (RLTp).}
	\label{fig:input:pdf_rltp}
\end{figure}

\paragraph{Retrolaminar tissue pressure.}
For the RLTp we used a normal distribution (Fig.~\ref{fig:input:pdf_rltp}) with mean $\mu_{normal}=9.5 \, \unit{mmHg}$ and standard deviation $\sigma_{normal}=2.2\, \unit{mmHg}$~\cite{ren2010cerebrospinal}.
In this case, no issues regarding the independence with other input or the sampling of unphysiological values occur.

\begin{remark}
	Our choices for the input probability density functions are dictated by clinical measurements reported in the literature, in particular \cite{sesso2000systolic} for the normal distribution of MAP, \cite{suh2012distribution} for the lognormal distribution of IOP, and \cite{ren2010cerebrospinal} for the normal distribution of RLTp.
	For the blood pressure and the retrolaminar tissue pressure, all the clinical literature refers to a normal probability density function of these inputs, whereas for what concerns the IOP, the prior knowledge on the uncertainty distribution is still a matter of debate. 
    The rationale behind our specific choices is that we aim to avoid unphysiological input values that may lead to unrealistic discontinuities in the simulation results. It would also be possible in the future to estimate the input probability density functions in a more patient-specific manner, by using given repeated clinical measurements of patient-specific targets, such as systolic and diastolic blood pressure, intracranial pressure, intraocular pressure etc.
\end{remark}


\section{Results}
\label{sec:results}
In this section, we present two virtual studies using the OMVS described in Sec.~\ref{subsec:math_comp} and employing the two methods described in Sec.~\ref{subsec:uq_sa}.
The input data assumptions have been discussed in Sec.~\ref{subsec:input_data}. 
In particular, in the first study, we will use uncertainty propagation, while in the second one, we will run a sensitivity analysis on the model. 

\subsection{Uncertainty propagation study} \label{subsec:uq-study}
We completed three different sets of $N=10000$ evaluations of the model using only the IOP as a stochastic input distribution.
For the RTLp we fixed its value for the three sets at $9.5 \, \unit{mmHg}$  (mean value provided by ~\cite{ren2010cerebrospinal}). 
These three sets of evaluations differ from each other by the blood pressure value imposed; we selected three cases of clinical interest - {in the same spirit} as in~\cite{guidoboni2014intraocular}:
\begin{enumerate}
	\item \textit{baseline} subjects with a systolic/diastolic blood pressure of $120/80$ \unit{mmHg};
	\item \textit{low} blood pressure subjects with SP $=100$ \unit{mmHg} and DP $=70$ \unit{mmHg};
	\item \textit{high} blood pressure subjects with SP/DP $= 140/90$ \unit{mmHg}.
\end{enumerate}
Recall that the outputs on which we focus for the UQ study are listed in Tab.~\ref{tab:output_params}.

\paragraph{Numerical simulations.}

\begin{figure}[h!]
	\centering
	\begin{subfigure}{0.32\linewidth}
		\includegraphics[width=\linewidth]{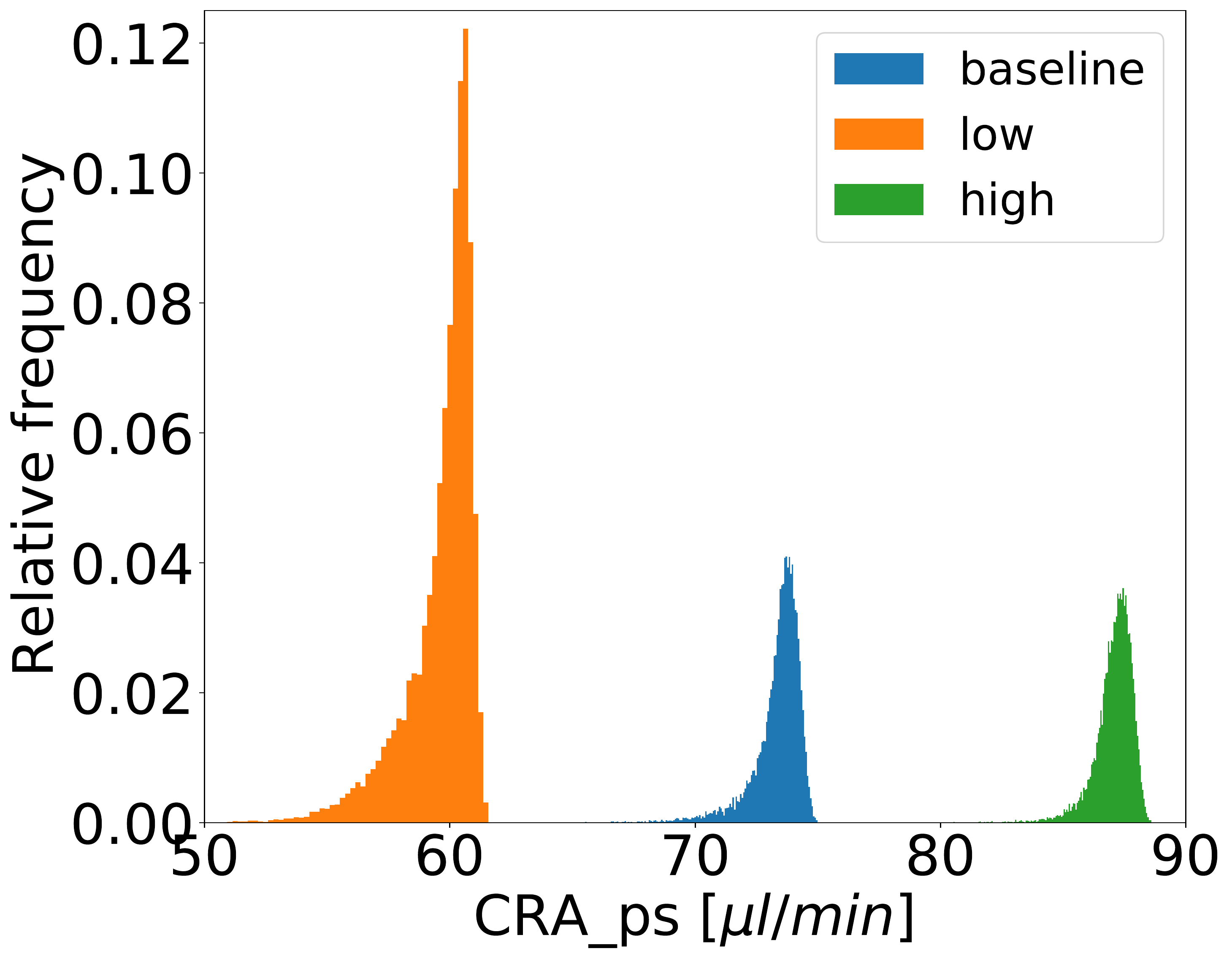}
		\caption{{Peak systolic CRA blood flow.}}
		\label{fig:uq:CRA:ps}
	\end{subfigure}
	\begin{subfigure}{0.32\linewidth}
		\includegraphics[width=\linewidth]{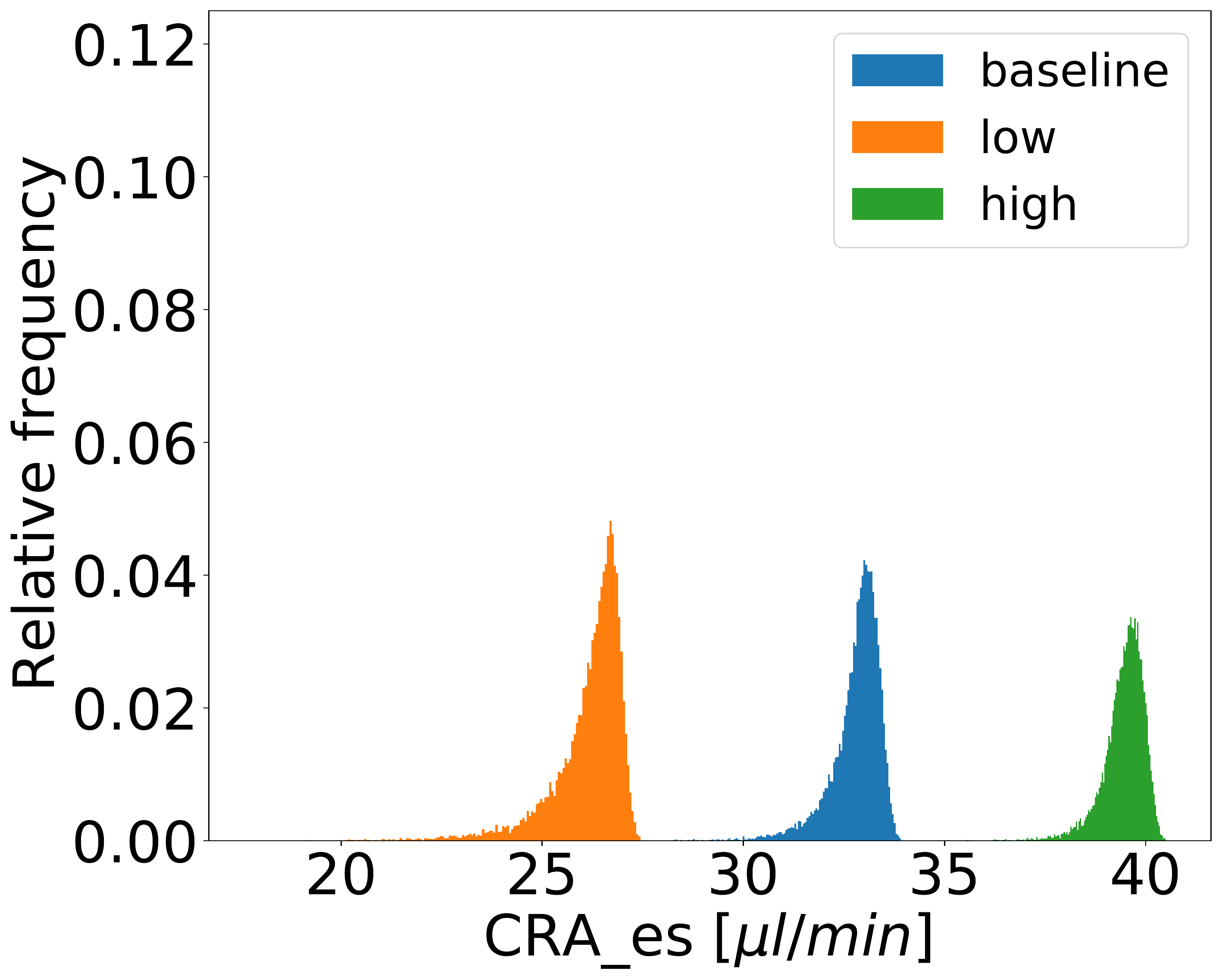}
		\caption{{End systolic CRA blood flow.}}
		\label{fig:uq:CRA:es}
	\end{subfigure}
	\begin{subfigure}{0.32\linewidth}
		\includegraphics[width=\linewidth]{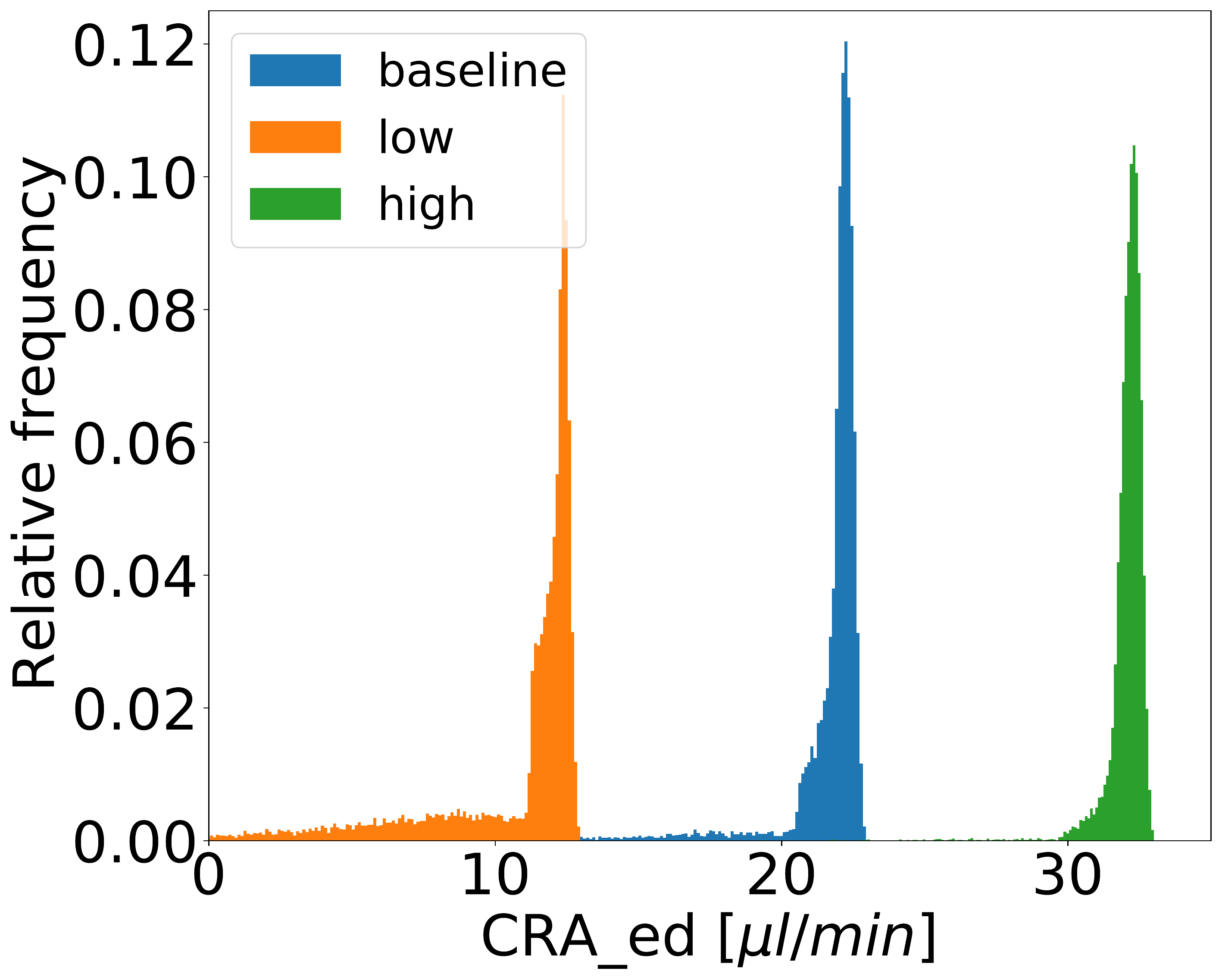}
		\caption{{End diastolic CRA blood flow.}}
		\label{fig:uq:CRA:ed}
	\end{subfigure}
	\begin{subfigure}{0.32\linewidth}
		\includegraphics[width=\linewidth]{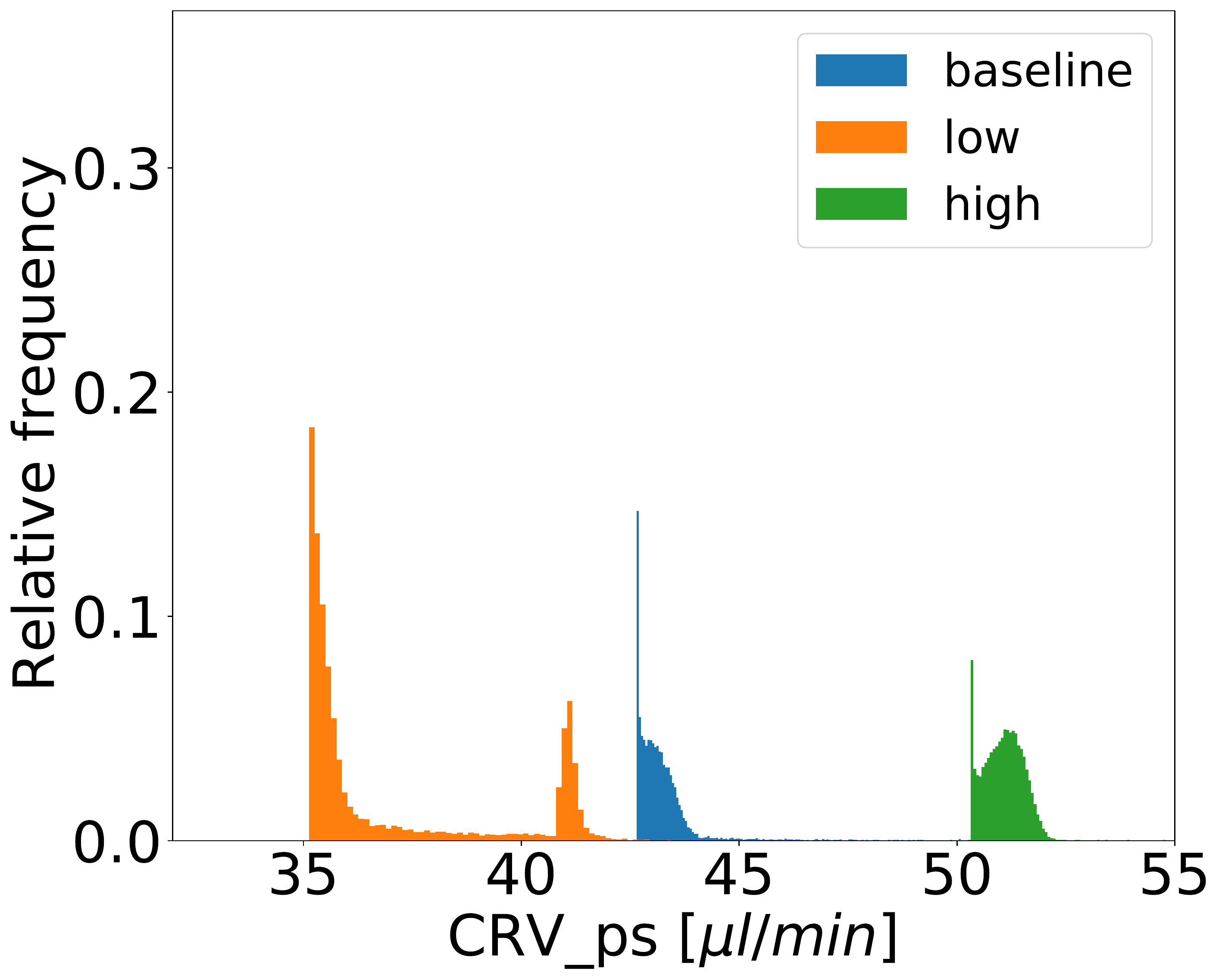}
		\caption{\review{Peak systolic CRV blood flow.}}
		\label{fig:uq:CRV:ps}
	\end{subfigure}
	\begin{subfigure}{0.32\linewidth}
		\includegraphics[width=\linewidth]{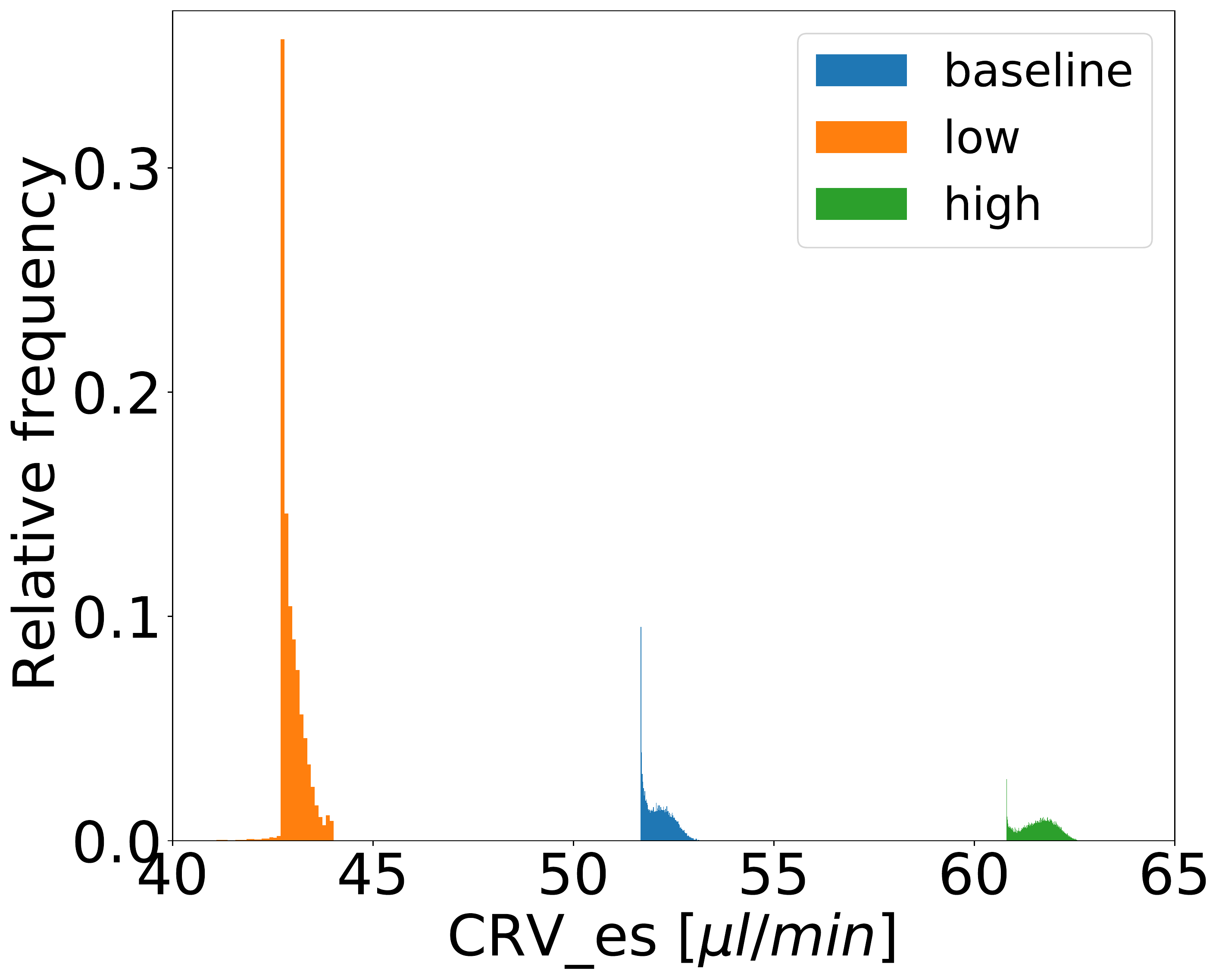}
		\caption{\review{End systolic CRV blood flow.}}
		\label{fig:uq:CRV:es}
	\end{subfigure}
	\begin{subfigure}{0.32\linewidth}
		\includegraphics[width=\linewidth]{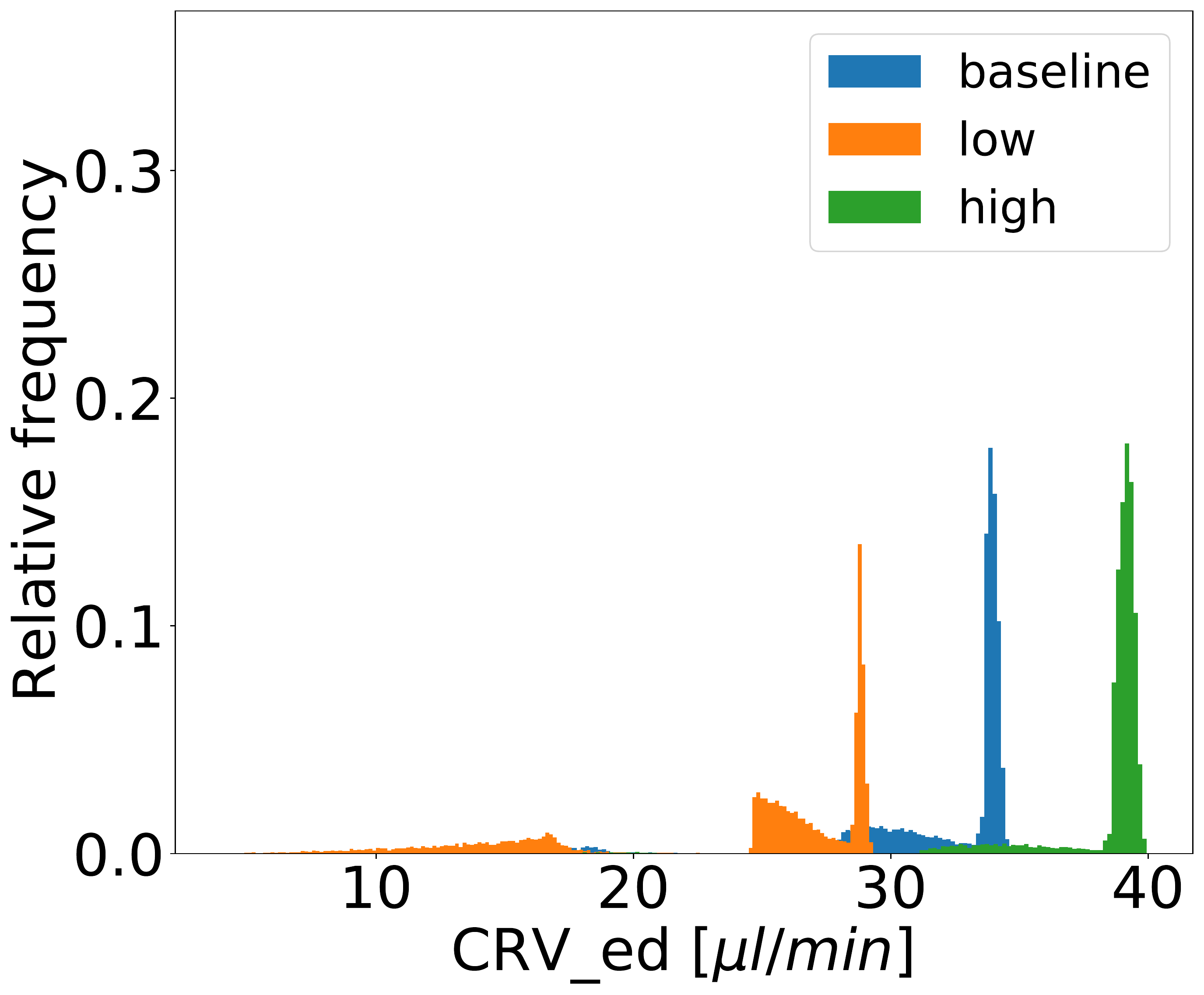}
		\caption{\review{End diastolic CRV blood flow.}}
		\label{fig:uq:CRV:ed}
	\end{subfigure}
		\begin{subfigure}{0.32\linewidth}
		\includegraphics[width=\linewidth]{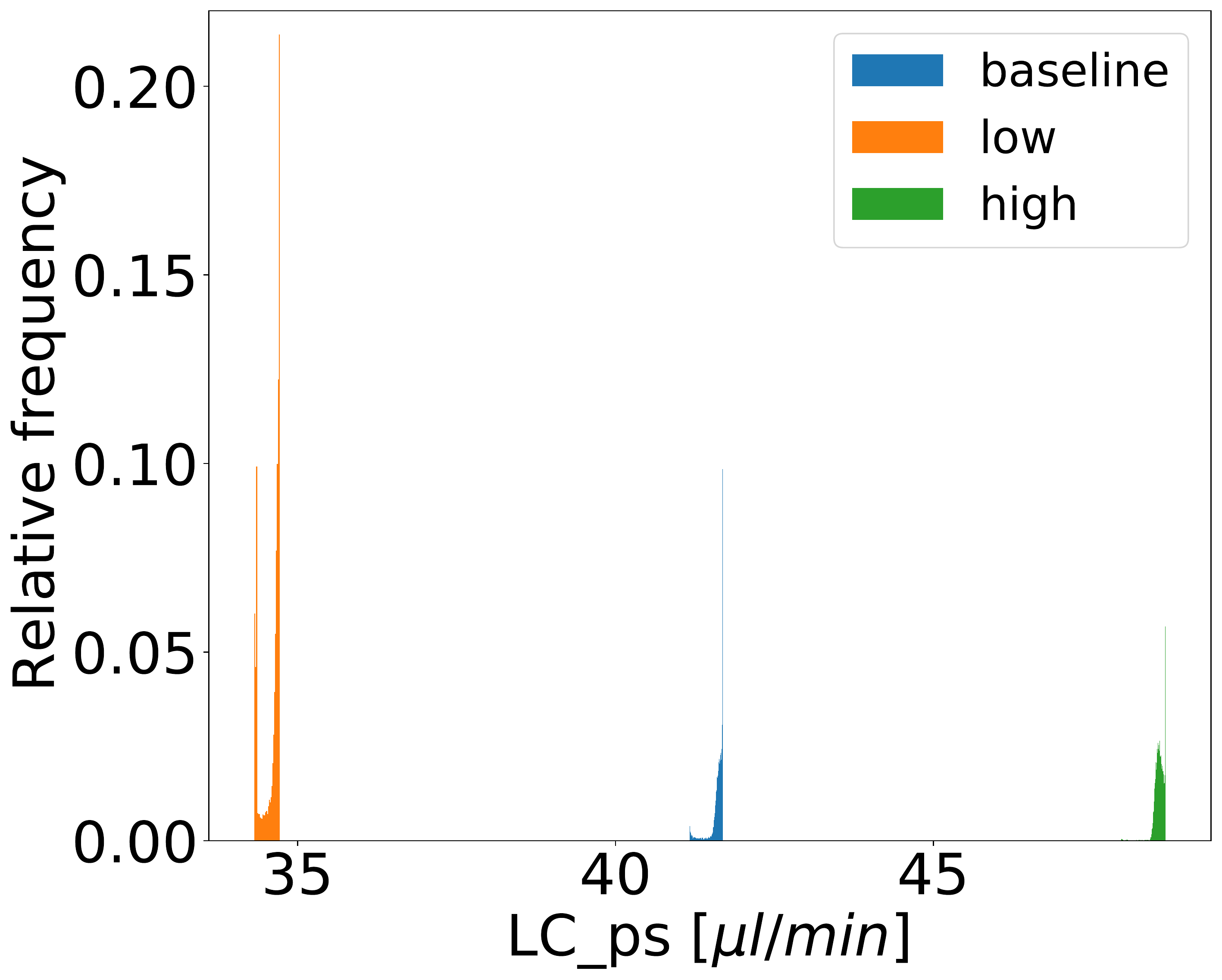}
		\caption{\review{Peak systolic LC blood flow.}}
		\label{fig:uq:LC:ps}
	\end{subfigure}
	\begin{subfigure}{0.32\linewidth}
		\includegraphics[width=\linewidth]{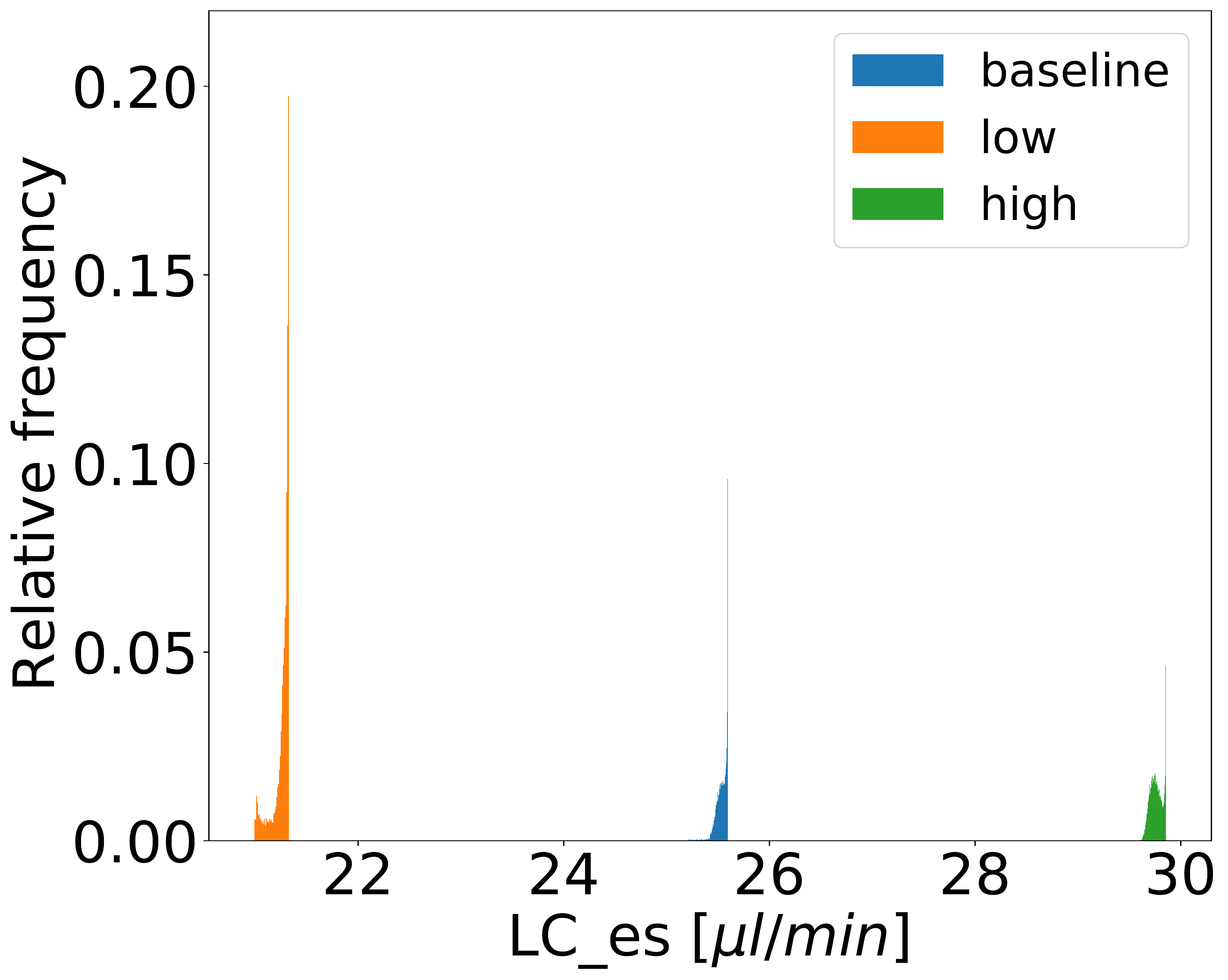}
		\caption{\review{End systolic LC blood flow.}}
		\label{fig:uq:LC:es}
	\end{subfigure}
	\begin{subfigure}{0.32\linewidth}
		\includegraphics[width=\linewidth]{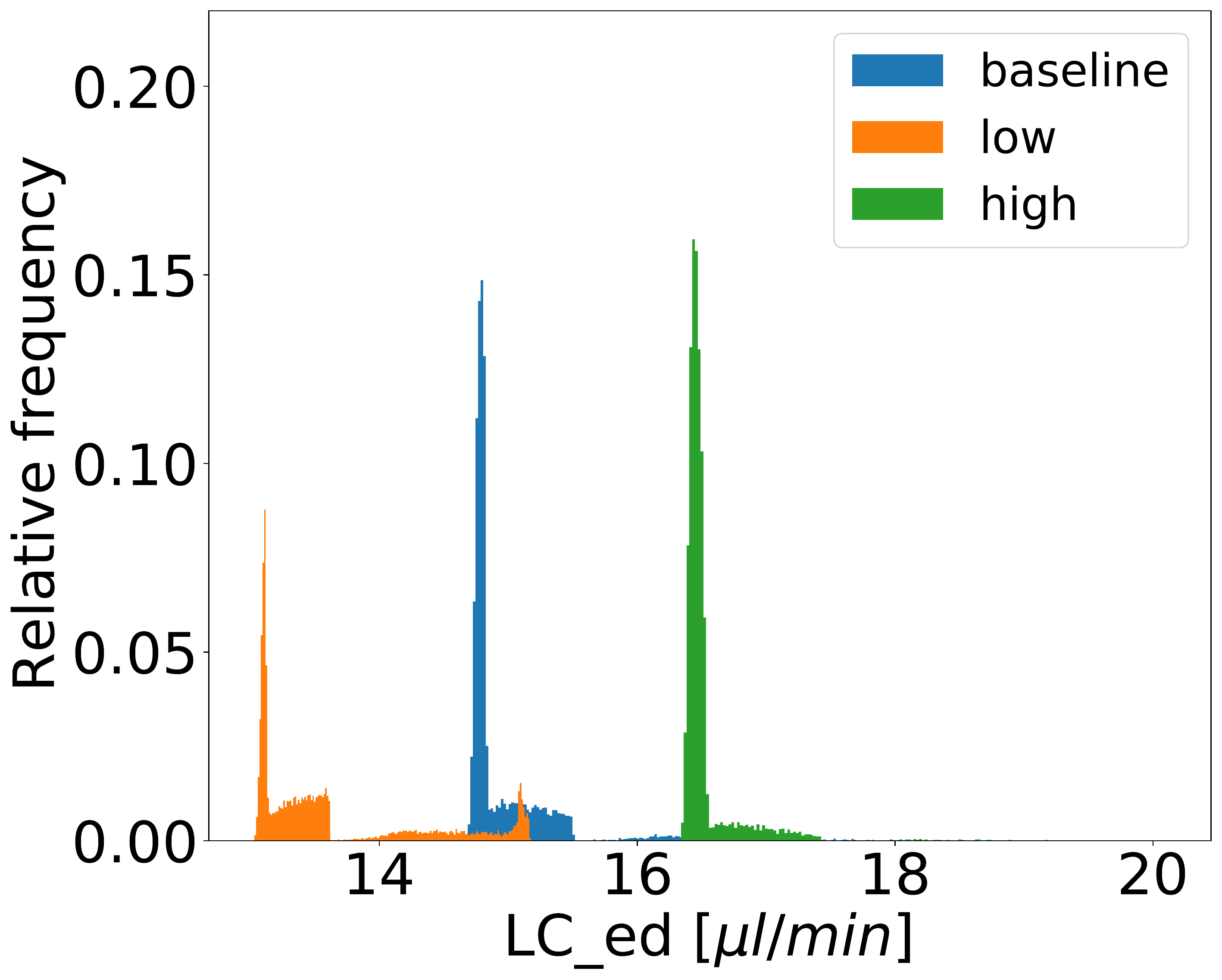}
		\caption{\review{End diastolic LC blood flow.}}
		\label{fig:uq:LC:ed}
	\end{subfigure}
	\caption{Output probability density functions \review{predicted by the OMVS for three virtual populations: baseline subjects with normal SP and DP values (blue), low blood pressure subjects (orange), and high blood pressure subjects (green).}}
	\label{fig:uq:pdf}
\end{figure}

\begin{figure}[h]
	\centering
	\begin{subfigure}{0.45\linewidth}
		\centering
		\includegraphics[width=\linewidth]{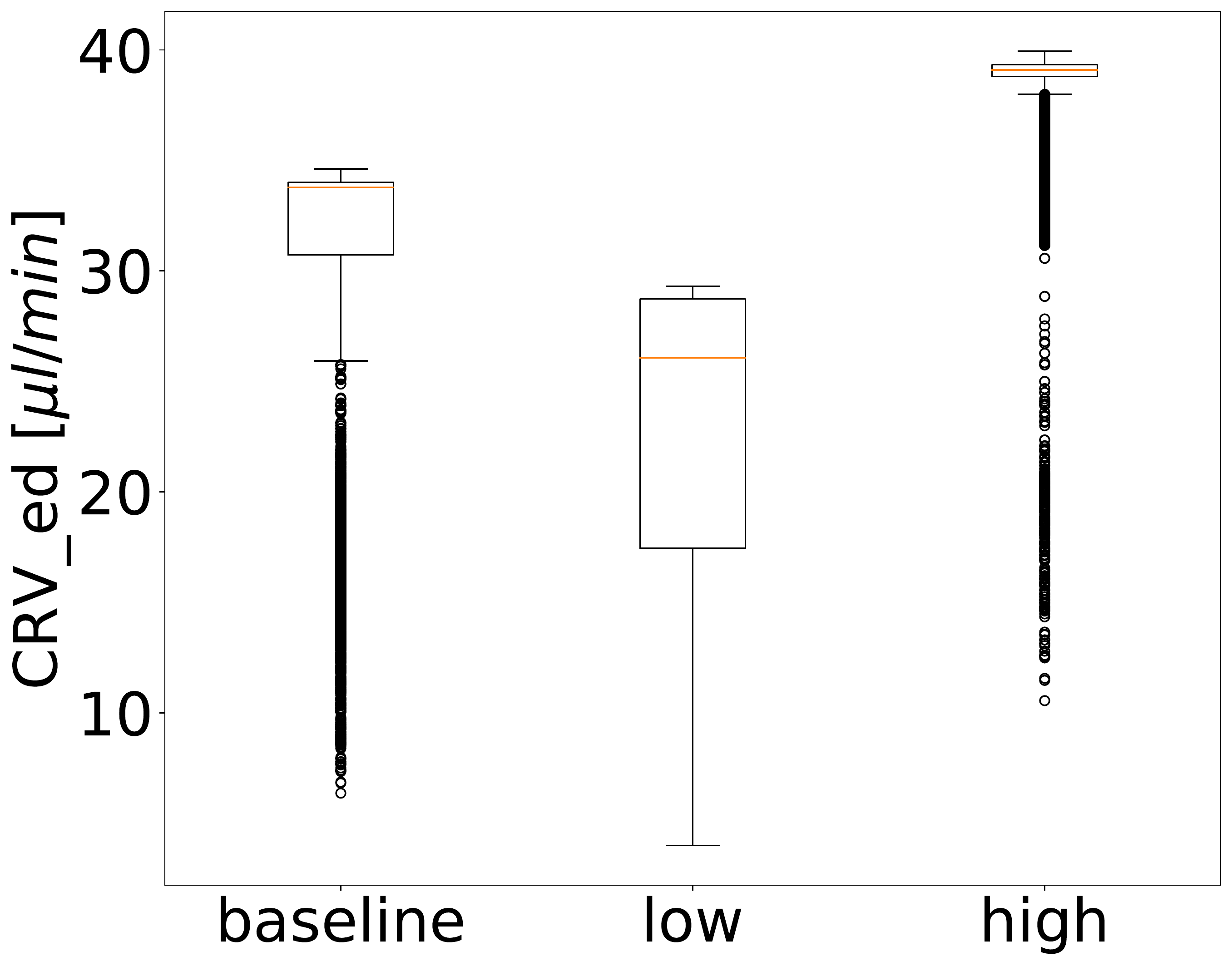}
		\caption{{End diastolic CRV blood flow.}}
		\label{fig:uq:CRVboxplot:ed}
	\end{subfigure}
	\begin{subfigure}{0.45\linewidth}
		\centering
		\includegraphics[width=\linewidth]{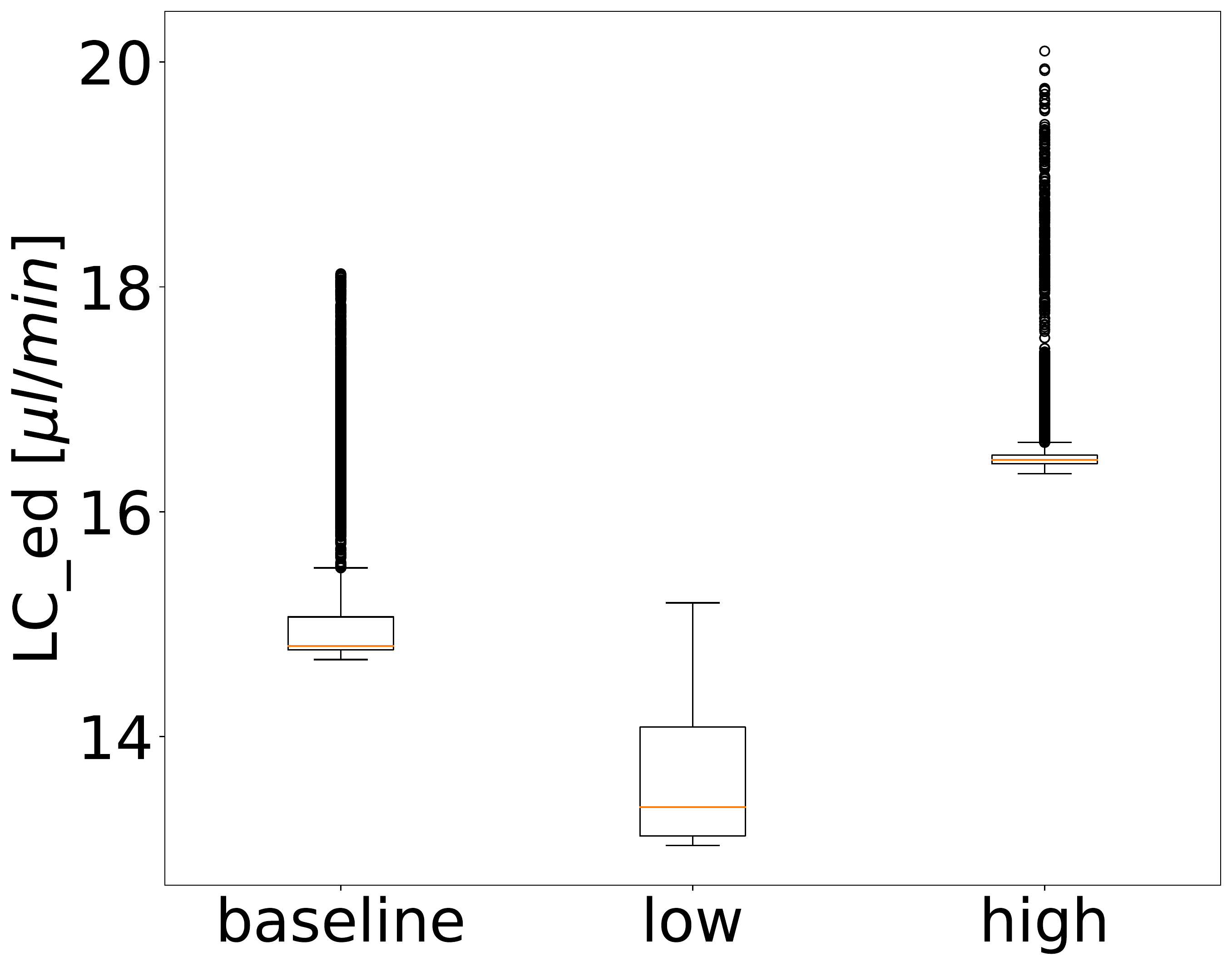}
		\caption{{End diastolic LC blood flow.}}
		\label{fig:uq:LCboxplot:ed}
	\end{subfigure}
	\caption{Boxplots of the end diastolic blood flow in the CRV and the lamina cribrosa predicted by the OMVS for three virtual populations: \textit{baseline} subjects with normal SP and DP values, \textit{low} blood pressure subjects, and \textit{high} blood pressure subjects.}
	\label{fig:boxplot}
\end{figure}

Fig.~\ref{fig:uq:pdf} shows the computed probability density functions for the three different locations we selected (CRA, CRV and LC), for three different time instants (peak systole, end systole and end diastole), and for the three populations of clinical interest (low, baseline and high blood pressure). We also report the simulated mean and standard deviation for each computed output in Tab.~\ref{tab:uq:results}.
Thus, we highlight that:
\begin{itemize}
	\item as expected, in each location and for all time instants, blood flow values decrease as we move from \textit{low} to \textit{baseline} and \textit{high} blood pressure populations;
	\item for the simulated CRA blood flow, the three populations (\textit{baseline}, \textit{low} and \textit{high}) have distinct pdfs at all evaluated time instants (peak systolic, end systolic, end diastolic), and present a significant asymmetry;
	\item for the simulated CRV blood flow, we have well established different pdfs at peak and end systole. 
	We remark that within the veins we have a delayed peak of blood flows and CRV\_es values higher than CRV\_ps values (Tab.~\ref{tab:uq:results}), in good agreement with experimentally observed patterns in time velocity curves acquired with Doppler imaging.
	The CRV blood flow pdf at end diastole exhibits a peculiar shape: a peak of frequencies for high values and a plateau in the number of realizations for low values separated by an almost empty frequency area of relatively middle values. 
	Even if the clinical interpretation of these distributions is difficult, it is interesting to observe that these results suggest a different repartition of frequencies between the three cases, in particular, for \textit{low} blood pressures, the high values peak is narrower and the plateau is wider and with more realizations than the \textit{high} blood pressure case. This fact is confirmed by the boxplot in Fig.~ \ref{fig:uq:CRVboxplot:ed}: the tail of CRV blood flow low values is within the first and third quantile range for the \textit{low} blood pressure case, whereas this tail is composed just by outliers for \textit{baseline} and \textit{high} blood pressure populations;
	\item for the LC blood flow, the OMVS suggests a similar analysis to CRV blood flow. The pdfs at peak, and end systole are distinct, whereas at end diastole it shows a peak of frequencies at low values and a more uniform distribution elsewhere. 
	Also, Fig.~\ref{fig:uq:LC:ed} points out that for only for the \textit{low} blood pressure virtual population a second peak and a considerable high number of frequencies can be identified at high values. 
	Similarly to the CRV analysis, we propose the boxplot for the output LC\_ed in Fig.~\ref{fig:uq:LCboxplot:ed}: as predicted by the pdf. The tail of high values in the \textit{low} blood pressure population is within the first and third quantile range. In contrast, the same tail consists of only outliers for \textit{baseline} and \textit{high} blood pressure populations.
	Finally, the computed LC blood flow  variability is considerably lower than for the other two outputs (CRA and CRV) as reported in Tab.~\ref{tab:uq:results}.
\end{itemize} 

\begin{table}[h!]
	\centering
	\resizebox{\linewidth}{!}{
	\begin{tabular}{l | c c c c c c c c c}
		$[\mu l/min]$& \textsc{CRA\_ps} & \textsc{CRA\_es} & \textsc{CRA\_ed} & \textsc{CRV\_ps} & \textsc{CRV\_es} & \textsc{CRV\_ed} & \textsc{LC\_ps} & \textsc{LC\_es} & \textsc{LC\_ed} \\ \hline
		\textit{baseline} & $73.3 \pm 1.0$ & $32.8 \pm 0.6$ & $21.7 \pm 1.7$ & $43.6 \pm 1.7$ & $52.1 \pm 0.3$ & $31.7 \pm 4.6$ & $41.6 \pm 0.1$ & $25.5 \pm 0.05$ & $15.0 \pm 0.5$ \\
		\textit{low} & $59.6 \pm 1.5$ & $26.1 \pm 0.9$ & $10.6 \pm 3.0$ & $37.0 \pm 2.4$ & $42.9 \pm 0.3$ & $23.7 \pm 6.3$ & $34.6 \pm 0.1$ & $21.3 \pm 0.08$ & $13.6 \pm 0.7$ \\
		\textit{high} & $87.1 \pm 0.8$ & $39.5 \pm 0.5$ & $32.1 \pm 0.8$ & $51.1 \pm 0.9$ & $61.6 \pm 0.4$ & $38.3 \pm 2.7$ & $48.6 \pm 0.07$ & $29.8 \pm 0.05$ & $16.5 \pm 0.3$ \\
	\end{tabular}}
	\caption{Simulated mean and standard deviation for all the quantities of interest in the uncertainty propagation study.}
	\label{tab:uq:results}
\end{table}


\subsection{Sobol' index study} \label{subsec:sobol-study}

In this study, we utilize as random input variables the IOP, the retrolaminar tissue pressure (RLTp), and the SP and DP to compute the Sobol' indices using the distribution introduced in Sec.~\ref{subsec:input_data}. 
We compute the Sobol' first and total indices using the Saltelli algorithm~\cite{saltelli2002} and the FAST first and total indices~\cite{saltelli1999}.
{We performed 5 analysis with increasing $\mathcal{N} = \left[ 1000, \: 2000, \: 5000, \: 7500, \:10000\right]$. 
The stopping criteria is based on the absolute iterative error:
\begin{equation}
 \max_{\forall o \in \mathcal{O}} \left\{ |idx_n - idx_{n^*}| \right\}
 \label{eq:error}
\end{equation}
where $\mathcal{O}= CRA_{ps,es,ed}, \,CRV_{ps,es,ed}, \, LC_{ps,es,ed}$, $idx$ are the first or total index for two consecutive choices ($n > n^*$) in $\mathcal{N}$.  
The final figures proposed have been obtained with $n = 10000$ where the error was less than $4 \cdot 10^{-2}$ for all input indices and using both algorithms to compute the results (see Tab. \ref{tab:sobol:err}).}

\begin{table}[h!]
	\centering
	\resizebox{0.9\linewidth}{!}{{
	\begin{tabular}{c p{1mm} c c p{1mm} c c}
		\textbf{Input} && \multicolumn{2}{c}{\textbf{Monte-Carlo approach~\cite{saltelli2002}}} && \multicolumn{2}{c}{\textbf{FAST method~\cite{saltelli1999}}} \\
		&& \textsc{first order index} & \textsc{total order index} && \textsc{first order index} & \textsc{total order index} \\[5pt]
		\multicolumn{7}{c}{$n=2000, \; n^* = 1000$} \\ \hline
		\textit{IOP} && 0.332 & 0.009 && 0.012 & 0.012 \\
		\textit{RLTp} && 0.056 & 0.007 && 0.0004 & 0.006 \\
		\textit{MAP} && 0.092 & 0.090 && 0.032 & 0.0192 \\[5pt]
		\multicolumn{7}{c}{$n=5000, \; n^* = 2000$} \\ \hline
		\textit{IOP} && 0.129 & 0.017 && 0.035 & 0.017 \\
		\textit{RLTp} && 0.058 & 0.006 && 0.0004 & 0.008 \\
		\textit{MAP} && 0.076 & 0.063 && 0.031 & 0.012 \\[5pt]
		\multicolumn{7}{c}{$n=7500, \; n^* = 5000$} \\ \hline
		\textit{IOP} && 0.038  & 0.029 && 0.034 & 0.016 \\
		\textit{RLTp} && 0.028 & 0.003 && 0.0003 & 0.013 \\
		\textit{MAP} && 0.046 & 0.039 && 0.004 & 0.004 \\[5pt]		
		\multicolumn{7}{c}{$n=10000, \; n^* = 7500$} \\ \hline
		\textit{IOP} && 0.037  & 0.014 && 0.006 & 0.005 \\
		\textit{RLTp} && 0.025 & 0.01 && 0.0002 & 0.006 \\
		\textit{MAP} && 0.02 & 0.017 && 0.007 & 0.007 \\[5pt]		
	\end{tabular}}}
	\caption{\review{Error based on Eq.~(\ref{eq:error}) computed between two consecutive choices of $n \in \mathcal{N}$.}}
	\label{tab:sobol:err}
\end{table}

\begin{figure}[h!]
	\centering
	\begin{subfigure}{0.32\linewidth}
		\includegraphics[width=\linewidth]{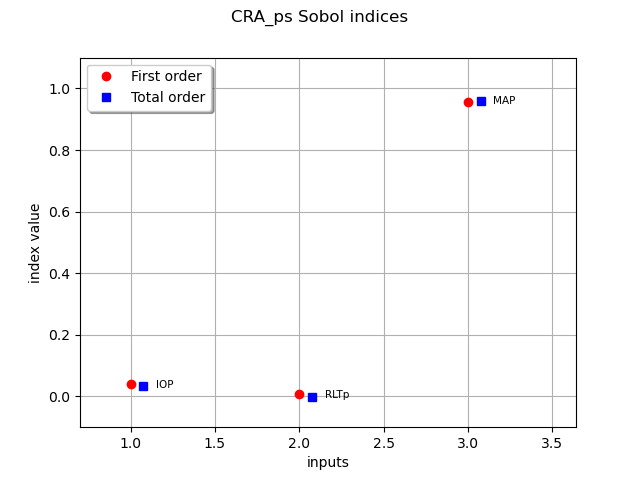}
		\caption{\review{Peak systolic CRA blood flow.}}
		\label{fig:sobol:CRA:ps}
	\end{subfigure}
	\begin{subfigure}{0.32\linewidth}
		\includegraphics[width=\linewidth]{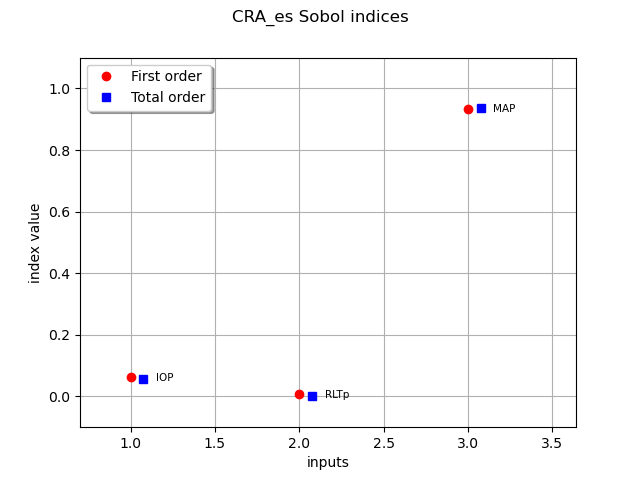}
		\caption{\review{End systolic CRA blood flow.}}
		\label{fig:sobol:CRA:es}
	\end{subfigure}
	\begin{subfigure}{0.32\linewidth}
		\includegraphics[width=\linewidth]{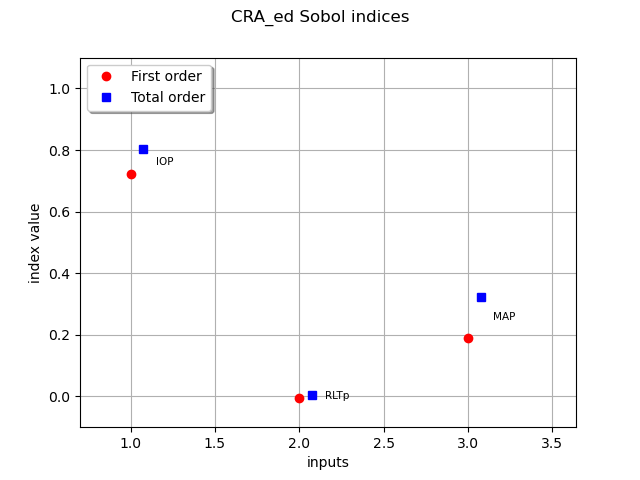}
		\caption{\review{End diastolic CRA blood flow.}}
		\label{fig:sobol:CRA:ed}
	\end{subfigure}
	\begin{subfigure}{0.32\linewidth}
		\includegraphics[width=\linewidth]{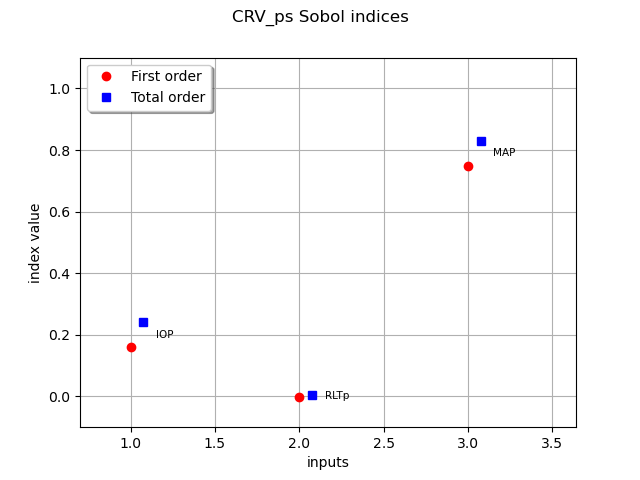}
		\caption{\review{Peak systolic CRV blood flow.}}
		\label{fig:sobol:CRV:ps}
	\end{subfigure}
	\begin{subfigure}{0.32\linewidth}
		\includegraphics[width=\linewidth]{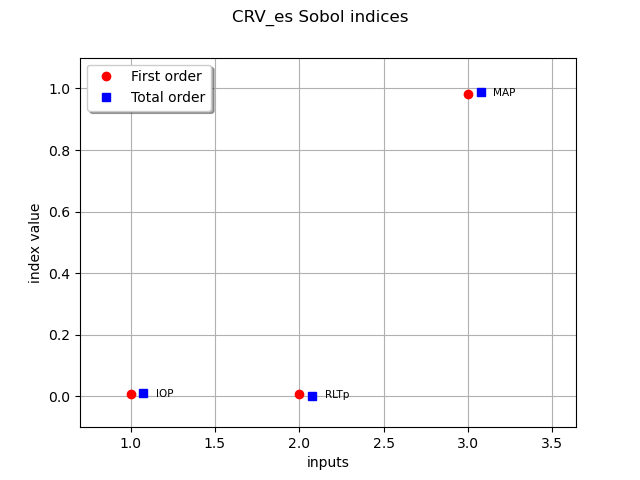}
		\caption{\review{End systolic CRV blood flow.}}
		\label{fig:sobol:CRV:es}
	\end{subfigure}
	\begin{subfigure}{0.32\linewidth}
		\includegraphics[width=\linewidth]{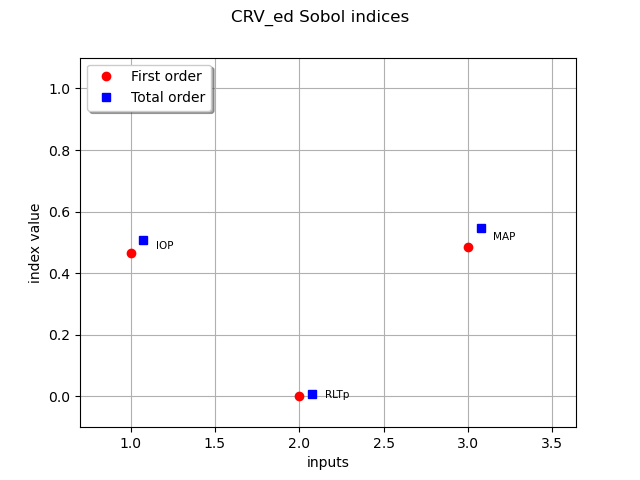}
		\caption{\review{End diastolic CRV blood flow.}}
		\label{fig:sobol:CRV:ed}
	\end{subfigure}
	\begin{subfigure}{0.32\linewidth}
		\includegraphics[width=\linewidth]{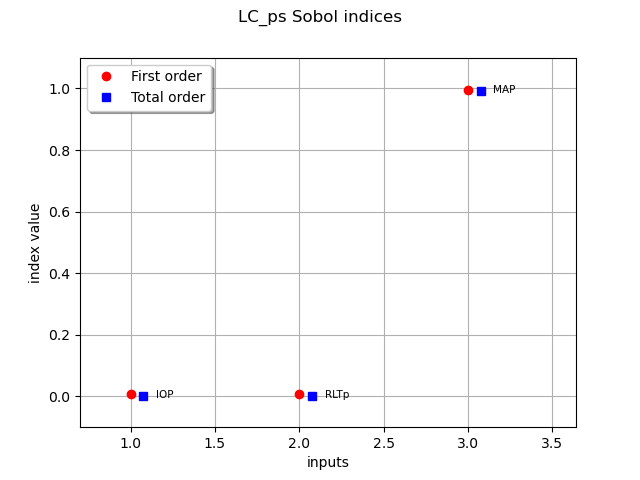}
		\caption{\review{Peak systolic LC blood flow.}}
		\label{fig:sobol:LC:ps}
	\end{subfigure}
	\begin{subfigure}{0.32\linewidth}
		\includegraphics[width=\linewidth]{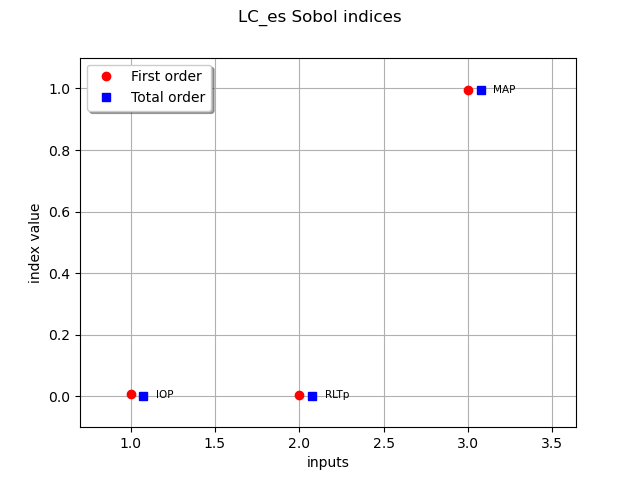}
		\caption{\review{End systolic LC blood flow.}}
		\label{fig:sobol:LC:es}
	\end{subfigure}
	\begin{subfigure}{0.32\linewidth}
		\includegraphics[width=\linewidth]{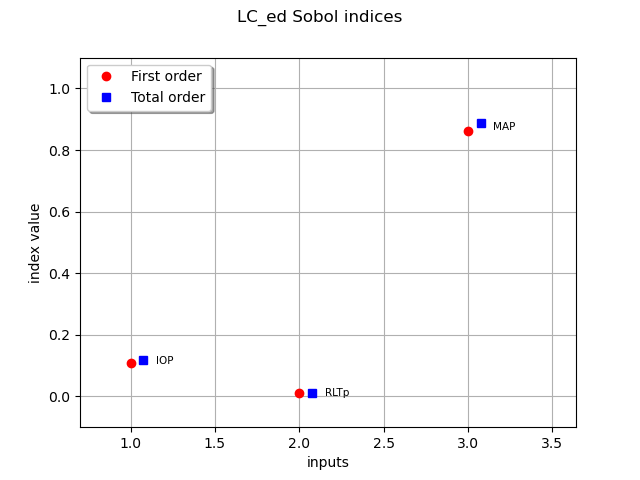}
		\caption{\review{End diastolic LC blood flow.}}
		\label{fig:sobol:LC:ed}
	\end{subfigure}
	\caption{Estimation of Sobol' indices using the Monte-Carlo approach~\cite{saltelli2002}. Red squares represent first order indices, while blue squares represent total order indices.}
	\label{fig:sobol}
\end{figure}

\begin{figure}[h!]
	\centering
	\begin{subfigure}{0.32\linewidth}
		\includegraphics[width=\linewidth]{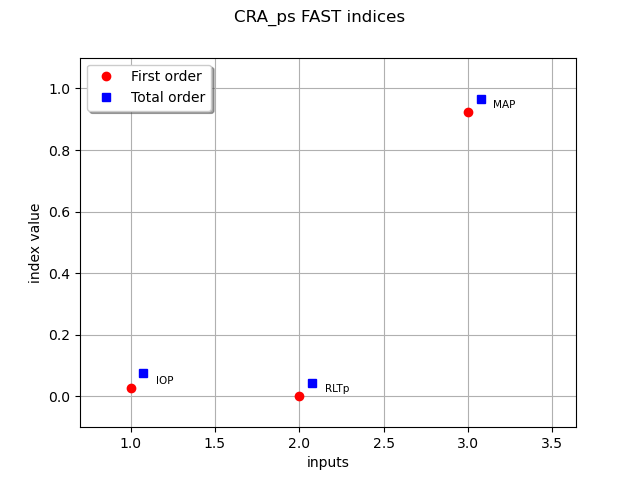}
		\caption{\review{Peak systolic CRA blood flow.}}
		\label{fig:fast:CRA:ps}
	\end{subfigure}
	\begin{subfigure}{0.32\linewidth}
		\includegraphics[width=\linewidth]{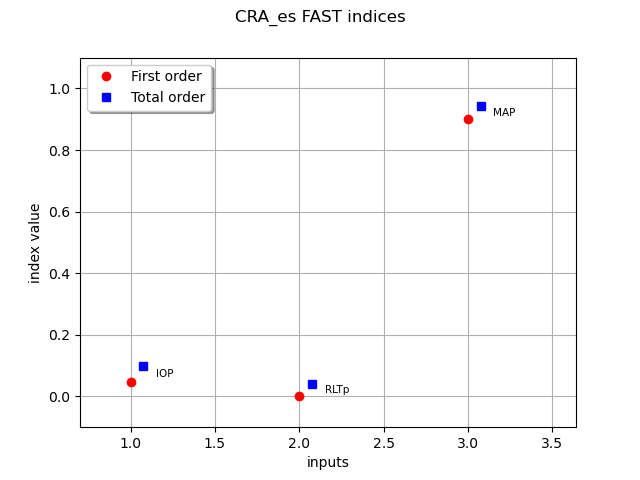}
		\caption{\review{End systolic CRA blood flow.}}
		\label{fig:fast:CRA:es}
	\end{subfigure}
	\begin{subfigure}{0.32\linewidth}
		\includegraphics[width=\linewidth]{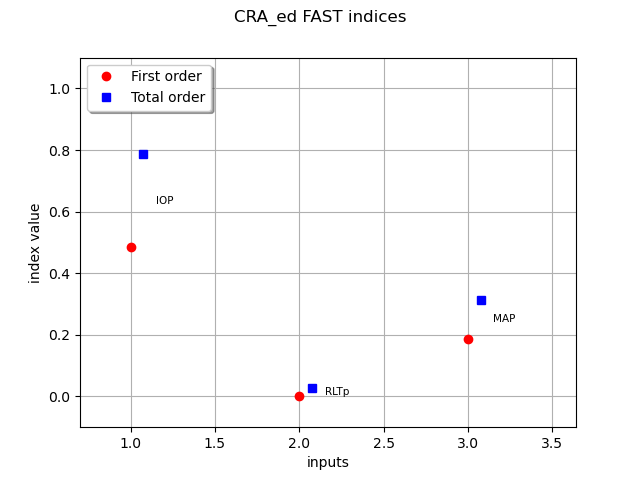}
		\caption{\review{End diastolic CRA blood flow.}}
		\label{fig:fast:CRA:ed}
	\end{subfigure}
	\begin{subfigure}{0.32\linewidth}
		\includegraphics[width=\linewidth]{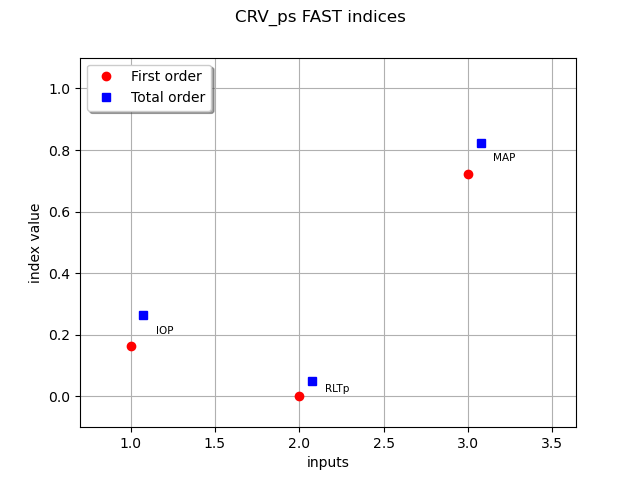}
		\caption{\review{Peak systolic CRV blood flow.}}
		\label{fig:fast:CRV:ps}
	\end{subfigure}
	\begin{subfigure}{0.32\linewidth}
		\includegraphics[width=\linewidth]{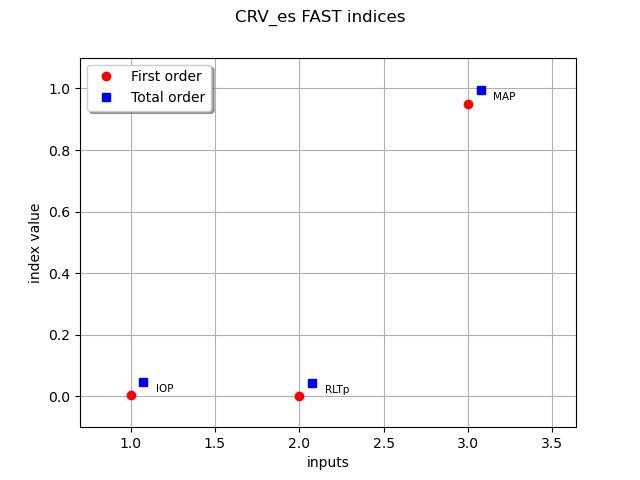}
		\caption{\review{End systolic CRV blood flow.}}
		\label{fig:fast:CRV:es}
	\end{subfigure}
	\begin{subfigure}{0.32\linewidth}
		\includegraphics[width=\linewidth]{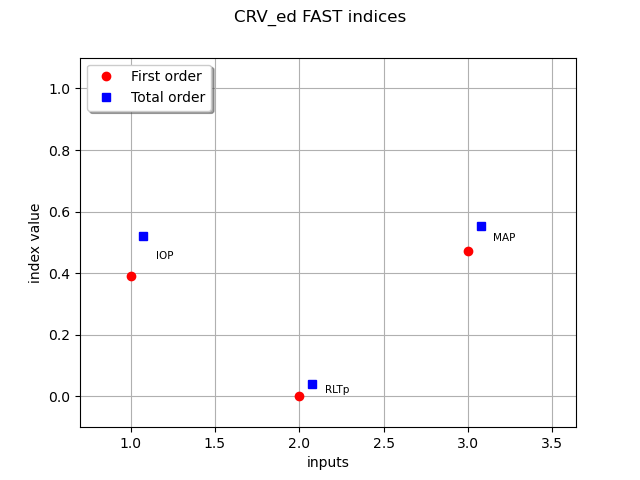}
		\caption{\review{End diastolic CRV blood flow.}}
		\label{fig:fast:CRV:ed}
	\end{subfigure}
	\begin{subfigure}{0.32\linewidth}
		\includegraphics[width=\linewidth]{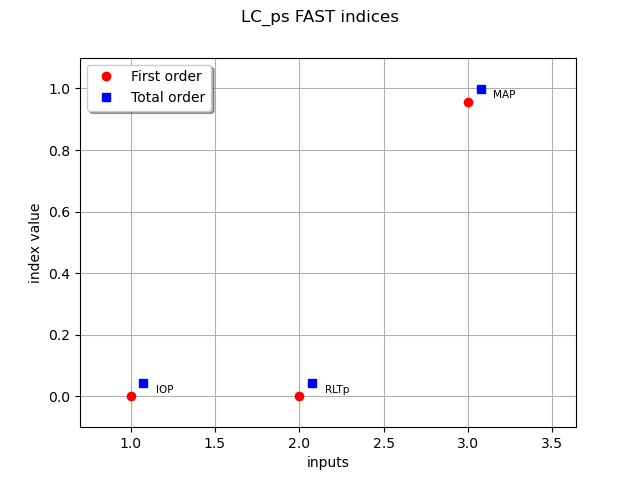}
		\caption{{Peak systolic LC blood flow.}}
		\label{fig:fast:LC:ps}
	\end{subfigure}
	\begin{subfigure}{0.32\linewidth}
		\includegraphics[width=\linewidth]{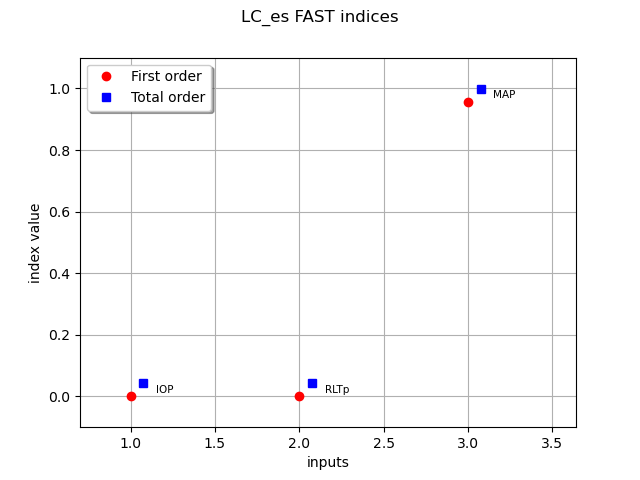}
		\caption{{End systolic LC blood flow.}}
		\label{fig:fast:LC:es}
	\end{subfigure}
	\begin{subfigure}{0.32\linewidth}
		\includegraphics[width=\linewidth]{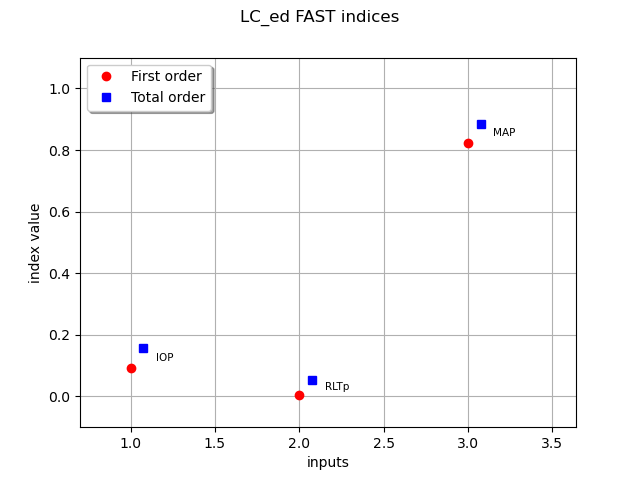}
		\caption{{End diastolic CRV blood flow.}}
		\label{fig:fast:LC:ed}
	\end{subfigure}
	\caption{Estimation of first and total indices using the FAST method~\cite{saltelli1999}. Red squares represent first order indices, while blue squares represent total order indices.}
	\label{fig:fast}
\end{figure}

\paragraph{Numerical simulations.}
Fig.~\ref{fig:sobol} report the Sobol indices using the Saltelli algorithm, while Fig.~\ref{fig:fast} report the FAST indices.

For what concerns the Sobol' indices, Figs.~\ref{fig:sobol:CRA:ps} and ~\ref{fig:sobol:CRA:es} point out that the CRA blood flow at peak and end systole is highly dependent on the value of the MAP, whereas  the influence of IOP and RTLp remains minimal.
In Fig.~\ref{fig:sobol:CRA:ed}, the results provided by the OMVS suggest that the CRA blood flow at end diastole is highly influenced by the IOP, moderately by the MAP, and almost negligibly by the RLTp. 
Moreover, in this case, we notice that the total order and first order, especially for the MAP, are significantly different, meaning that high order interactions among parameters contribute to the variance of this output. \\
For the CRV blood flow, Fig.~\ref{fig:sobol:CRV:ps} illustrates that there is a high dependency of CRV\_ps on MAP and only a moderate one from IOP.
CRV\_es depends mainly on MAP (Fig.~\ref{fig:sobol:CRV:es}), and CRV\_ed is highly dependent on the IOP and only mildly on  the other two inputs (Fig.~\ref{fig:sobol:CRV:ed}).
Also, for CRV blood flow at end diastole, we notice that high order interactions  occur, especially for MAP and IOP. \\
For the lamina cribrosa blood flow, the results provided by the OMVS sensitivity analysis suggest that the MAP is the dominant factor with a moderate influence of the IOP only at end diastole (Fig.~\ref{fig:sobol:LC:ps}, ~\ref{fig:sobol:LC:es}, and ~\ref{fig:sobol:LC:ed}). \\

The FAST indices (Figs.~\ref{fig:fast}) suggest that CRA\_ps (Fig.~\ref{fig:fast:CRA:ps}), CRA\_es (Fig.~\ref{fig:fast:CRA:es}) and LC\_ed (Fig. ~\ref{fig:fast:LC:ed}) depend mainly on MAP and mildly on IOP and RLTp. \\
For what concerns the CRA\_ed (Fig.~\ref{fig:fast:CRA:ed}) and CRV\_ed (Fig.~\ref{fig:fast:CRV:ed}), they are highly influenced by the IOP, moderately by MAP and negligibly little by the RLTp.
In this case the first and total order - especially for the IOP - are very different, implying that high order interactions are quite relevant to explain the variability of the output. \\
Fig.~\ref{fig:fast:CRV:ps} suggests that CRV\_ps has a high dependency on MAP and a moderate one on IOP. \\
For CRV\_es (Fig.~\ref{fig:fast:CRV:es}), LC\_ps (Fig.~\ref{fig:fast:LC:ps}) and LC\_es (Fig. ~\ref{fig:fast:LC:es}), the sensitivity analysis results show that their variability is almost solely due to changes in MAP.

\section{Discussion.} \label{sec:discussion}
Hemodynamics in the ocular posterior tissue vasculature results from the combined effects of different factors. Specifically, ocular blood flow is driven by the difference between arterial and venous blood pressure, is impeded by IOP and RTLp (directly related to cerebrospinal fluid pressure) and is modulated by vascular regulation. Understanding these complex interrelated effects represents a major challenge when interpreting results from various clinical studies in ophthalmology. 
The current contribution introduces a network-based model in the framework of the Ocular Mathematical Simulator, that couples retinal blood circulation with a simplified description of the LC hemodynamics. 
Further on, the model is employed to theoretically assess the relative contribution of IOP, RTLp and MAP stochastic variations on several clinically meaningful outputs characterizing hemodynamics in the CRA, CRV and LC, by means of uncertainty propagation methods and variance-based sensitivity indexes.

The results suggest that the hemodynamic response of the CRA, CRV and LC vasculature to variations in IOP and RTLp presents noticeable differences among individuals with different blood pressures, that strongly influences all the computed outputs.
These model predictions are in good agreement with the experimental findings in~\cite{tham2018inter} and several theoretical studies~\cite{guidoboni2014intraocular,guidoboni2014effect,carichino2014}, designed to elucidate how and to which extent blood pressure is influencing the distribution of ocular hemodynamics. Note that in a clinical setting, there is an intrinsic difficulty of evaluating the contribution of each factor and the complex relationships among them.

The use of a lognormal distribution for the IOP translates into a non trivial interpretation of the output, especially at the level of the CRV and LC due to the non-linear character of the model and the complex interplay between factors (as pointed out in Sec. \ref{subsec:sobol-study}).
In this context, mathematical models are crucial to reproduce these mechanisms and help to unveil their interpretation.

We have compared the CRA blood flow simulated results to other data in the literature.
Our results - Tab. \ref{tab:uq:results}, in particular \textit{baseline} where CRA\_ps $=73.3 \pm 1.0\, \mu l/min$, CRA\_es $=32.8 \pm 0.6\, \mu l/min$, CRA\_ed $=21.7 \pm 1.7\, \mu l/min$ - are in agreement with the experimental results of Dorner et al. \cite{dorner2002calculation}  (CRA mean blood flow of $38.1 \pm 9.1 \, \mu l/min$) and Riva et al. \cite{riva1979bidirectional,riva1985blood}(CRA mean blood flow of $33 \pm 9.6 \, \mu l/min$).

To compare our results with other clinical and mathematical studies that are more focused on the CRA blood velocities, we set the hypothesis of CRA diameter of about $160 \, \mu m$ \cite{dorner2002calculation,lee2007association}.
Using this assumption, our simulations provide similar values for CRA blood flows than the one measured by Harris et al. \cite{harris1996acute} (CRA\_ps $=120.6\, \mu l/min$, CRA\_ed $=30.1\, \mu l/min$) and the three virtual populations simulated by Guidoboni et al. \cite{guidoboni2014intraocular} (baseline: CRA\_ps $=119.4\, \mu l/min$, CRA\_ed $=33.7\, \mu l/min$; low: CRA\_ps $=95.3\, \mu l/min$, CRA\_ed $=28.9\, \mu l/min$; high: CRA\_ps $=142.3\, \mu l/min$, CRA\_ed $=41.0\, \mu l/min$).
This comparison shows the quality of the results, albeit the simplicity of the model we have employed for our study.
From the SA (Sec. \ref{subsec:sobol-study}) we evince that, as expected, the CRA blood flow depends mainly by the MAP, and only at end diastole - when the arterial pressure is at minimum - the IOP is affecting the CRA results. \\

For the CRV, we compared the output CRV\_es results (Tab. \ref{tab:uq:results}) with total venous blood flow measurements reported in the literature.
The simulated baseline mean value ($=52.1 \, \mu l/min$) agrees with the experiments performed by Garcia et al. \cite{garcia2002retinal} ($64.9 \pm 12.8 \, \mu l/min$) and Feke et al. \cite{feke1978laser} ($80 \pm 12 \, \mu l/min$ for age 25-38 group, $73 \pm 13 \, \mu l/min$ for age 54-58 group).
From a qualitative viewpoint, the analysis of the low values tail proposed in the numerical results paragraph may have an interesting physiological interpretation.
As described in Section \ref{subsec:math_comp}, IOP has a non-linear effect on the retinal vasculature, particularly on the venous part.
Following this reason, the CRV high values peak may represent the natural state when the IOP is lower than the venous blood pressure, whereas the low values plateau denotes the collapse state.
This statement is consistent with the previous analysis where we found a uniform distribution for low blood pressures and an important frequency in the peak for high blood pressures. 
The SA supports this prediction: the dependency to IOP is significant not only at end diastole, but also at peak systole. 
In contrast with the CRA, the CRV is a venous vessel, which blood pressure is lower therefore more easily influenced by external pressures (\textit{e.g.} IOP).
These indices show also significant differences between the first-order and the total-order index, which means that there are high-order interactions between these parameters.
This fact is not surprising but it appears as a consequence of incorporating the Starling resistor effect,~\cite{starling1896absorption}, which is a crucial requirement to retrieve clinical data, as reviewed in~\cite{guidoboni2019}. \\

For the lamina cribrosa hemodynamics, we highlight the fact that all mathematical results are crucial in the investigation of disease because non-invasive measurements are not available nowadays for this tissue.
Following similar consideration made for the CRV, we notice an interesting physiological interpretation following the uncertainty propagation study (Sec. \ref{subsec:uq-study}).
For the \textit{low} blood pressure population, which is notably the most at risk for ocular neuropathies, the results may suggest that an overperfusion of the lamina -with respect to the mean values of that population -  occurs in more cases than for the other two population.
The simplicity of the model, as confirmed by the low variability in the LC blood flows, does not allow us to clarify this fact.
The SA (Sec. \ref{subsec:uq_sa}) supports this analysis.
Sobol and FAST indices show a high dependency only due to the MAP, which can be not so intuitive as for the CRA or the CRV, indeed we know that the pressure gradient across the LC ($IOP-RLTp$) may influence the hemodynamics \cite{sala2018analysis}.
Further analysis by adding a three-dimensional hemodynamical and biomechanical description of the lamina would certainly help in this investigation, in particular to notice the impact of IOP and RLTp on LC blood flow.\\

Finally we discuss the estimates of the  first and total indices using the Monte-Carlo and FAST approaches.
For all considered quantities of interest, the two indices are suggesting similar outcomes, which allows us also to cross-validate the results of our analysis.
The main differences concern CRA\_ed (Figs. \ref{fig:sobol:CRA:ed} and \ref{fig:fast:CRA:ed}) and CRV\_ed (Figs. \ref{fig:sobol:CRV:ed} and \ref{fig:fast:CRV:ed}), in particular on the high order interactions between IOP and MAP - more emphasized for the FAST estimates.
We privilege the results provided by the Monte-Carlo estimates which are unbiased compared to the FAST method. Also, the latter has the initial advantage of being more computationally efficient, but at the cost of extra assumptions of smoothness in the model \cite{xu2011reliability} that we are not currently able to verify for the OMVS.


\section{Conclusion}
\label{sec:concl}
Thanks to its special connection to the brain and its accessibility to measurements, the eye provides a unique window on the brain, thereby offering non-invasive access to a large set of potential biomarkers that might help in the early diagnosis and clinical care of Neuro-Degenerative Diseases \cite{guidoboni2020}.
The OMVS has already shown great potentiality to reproduce the ocular biomechanics and the hemodynamics \cite{sala2019mathematical}.

Pursuing this concept, in this contribution, we have proposed an uncertainty propagation study and a sensitivity analysis to evaluate the impact of uncertainties on this ophthalmological virtual laboratory.
First, we have set up a framework to perform a forward UQ analysis, which allowed us to evaluate the effects of the propagation of uncertainty from input to output. 
Second, we completed a SA study based on Sobol indices to capture the interplay between the different model parameters and their relative importance.
Finally, we have assessed qualitatively and quantitatively this computational framework in view of clinical applications.

In the context of computational models, a coupled UQ/SA analysis is crucial for the scientific research, especially in biology and medicine. 
The use of such mathematical framework may be employed in the process of product development or diseases understanding,  decreasing the number of physical tests necessary and therefore the economic cost. 
The low fidelity model employed already provided useful information for analysis, however this study must be pursued with higher fidelity models such as System II and System III, see Figs. ~\ref{fig:omvs_scheme:systemII} and ~\ref{fig:omvs_scheme:systemIII} respectively. 
	The current methodology could thus be further improved, in particular by 
\begin{inparaenum}[(\it i)]
	\item using a multilevel multifidelity estimator, as for instance the one developed in Dakota toolkit, see~ \cite{fleeter2020multilevel,adams2009dakota};
	\item devising a reduced order modelling approach for the 3D elastic and poroelastic models, following the reduced basis framework~\cite{Daversin2015,hild:tel-03025312}, developed in the open source software ~\texttt{Feel++} \cite{Prudhomme2012,Prudhomme2020}.
\end{inparaenum}
Finally, considering sensitivity indices over time, as in ~\cite{campos2020}, as well as characteristics instants, provides a promising perspective of the present work.


\small
\bibliography{biblio}

\providecommand{\href}[2]{#2}
\providecommand{\arxiv}[1]{\href{http://arxiv.org/abs/#1}{arXiv:#1}}
\providecommand{\url}[1]{\texttt{#1}}
\providecommand{\urlprefix}{URL }
\begin{thebibliography}{10}

\bibitem{baudin2017openturns}
\newblock M.~Baudin, A.~Dutfoy, B.~Iooss and A.-L. Popelin,
\newblock Openturns: An industrial software for uncertainty quantification in
  simulation,
\newblock \emph{Handbook of uncertainty quantification}, 2001--2038.

\bibitem{brault2017}
\newblock A.~Brault, L.~Dumas and D.~Lucor,
\newblock Uncertainty quantification of inflow boundary condition and proximal
  arterial stiffness--coupled effect on pulse wave propagation in a vascular
  network,
\newblock \emph{International journal for numerical methods in biomedical
  engineering}, \textbf{33} (2017), e2859.

\bibitem{campos2020}
\newblock J.~Campos, J.~Sundnes, R.~Dos~Santos and B.~Rocha,
\newblock Uncertainty quantification and sensitivity analysis of left
  ventricular function during the full cardiac cycle,
\newblock \emph{Philosophical Transactions of the Royal Society A},
  \textbf{378} (2020), 20190381.

\bibitem{carichino2014}
\newblock L.~Carichino, G.~Guidoboni, B.~Siesky, A.~Amireskandari,
  I.~Januleviciene, A.~Harris and P.~Causin,
\newblock Effect of intraocular pressure and cerebrospinal fluid pressure on
  the blood flow in the central retinal vessels,
\newblock \emph{Integrated Multidisciplinary Approaches in the Study and Care
  of the Human Eye Kugler Publications}, 59--66.

\bibitem{cassani2016blood}
\newblock S.~Cassani,
\newblock \emph{Blood circulation and aqueous humor flow in the eye:
  multi-scale modeling and clinical applications},
\newblock PhD thesis, 2016.

\bibitem{chen2013}
\newblock P.~Chen, A.~Quarteroni and G.~Rozza,
\newblock Simulation-based uncertainty quantification of human arterial network
  hemodynamics,
\newblock \emph{International journal for numerical methods in biomedical
  engineering}, \textbf{29} (2013), 698--721.

\bibitem{dorner2002calculation}
\newblock G.~T. Dorner, E.~Polska, G.~Garh{\"o}fer, C.~Zawinka, B.~Frank and
  L.~Schmetterer,
\newblock Calculation of the diameter of the central retinal artery from
  noninvasive measurements in humans,
\newblock \emph{Current eye research}, \textbf{25} (2002), 341--345.

\bibitem{Eck2015}
\newblock V.~Eck, J.~Feinberg, H.~Langtangen and L.~Hellevik,
\newblock Stochastic sensitivity analysis for timing and amplitude of pressure
  waves in the arterial system,
\newblock \emph{International journal for numerical methods in biomedical
  engineering}, \textbf{31} (2015), e02711.

\bibitem{feke1978laser}
\newblock G.~T. Feke and C.~E. Riva,
\newblock Laser doppler measurements of blood velocity in human retinal
  vessels,
\newblock \emph{JOSA}, \textbf{68} (1978), 526--531.

\bibitem{formaggia2010cardiovascular}
\newblock L.~Formaggia, A.~Quarteroni and A.~Veneziani,
\newblock \emph{Cardiovascular Mathematics: Modeling and simulation of the
  circulatory system}, vol.~1,
\newblock Springer Science \& Business Media, 2010.

\bibitem{fritzson2006openmodelica}
\newblock P.~Fritzson, P.~Aronsson, A.~Pop, H.~Lundvall, K.~Nystrom,
  L.~Saldamli, D.~Broman and A.~Sandholm,
\newblock Openmodelica-a free open-source environment for system modeling,
  simulation, and teaching,
\newblock in \emph{2006 IEEE Conference on Computer Aided Control System
  Design, 2006 IEEE International Conference on Control Applications, 2006 IEEE
  International Symposium on Intelligent Control},
\newblock IEEE, 2006,
\newblock 1588--1595.

\bibitem{garcia2002retinal}
\newblock J.~P. Garcia~Jr, P.~T. Garcia and R.~B. Rosen,
\newblock Retinal blood flow in the normal human eye using the canon laser
  blood flowmeter,
\newblock \emph{Ophthalmic research}, \textbf{34} (2002), 295--299.

\bibitem{gavish2008linear}
\newblock B.~Gavish, I.~Z. Ben-Dov and M.~Bursztyn,
\newblock Linear relationship between systolic and diastolic blood pressure
  monitored over 24 h: assessment and correlates,
\newblock \emph{Journal of hypertension}, \textbf{26} (2008), 199--209.

\bibitem{guidoboni2019}
\newblock G.~Guidoboni, A.~Harris and R.~Sacco,
\newblock Mathematical modeling of ocular fluid dynamics: From theory to
  clinical applications. modeling and simulation in science, engineering, and
  technology, 2019.

\bibitem{guidoboni2014effect}
\newblock G.~Guidoboni, A.~Harris, L.~Carichino, Y.~Arieli and B.~A. Siesky,
\newblock Effect of intraocular pressure on the hemodynamics of the central
  retinal artery: a mathematical model,
\newblock \emph{Mathematical Biosciences \& Engineering}, \textbf{11} (2014),
  523--546.

\bibitem{guidoboni2014intraocular}
\newblock G.~Guidoboni, A.~Harris, S.~Cassani, J.~Arciero, B.~Siesky,
  A.~Amireskandari, L.~Tobe, P.~Egan, I.~Januleviciene and J.~Park,
\newblock Intraocular pressure, blood pressure, and retinal blood flow
  autoregulation: a mathematical model to clarify their relationship and
  clinical relevance,
\newblock \emph{Investigative ophthalmology \& visual science}, \textbf{55}
  (2014), 4105--4118.

\bibitem{guidoboni2020}
\newblock G.~Guidoboni, R.~Sacco, M.~Szopos, L.~Sala, A.~C.
  Verticchio-Vercellin, B.~Siesky and A.~Harris,
\newblock Neurodegenerative disorders of the eye and of the brain: a
  perspective on their fluid-dynamical connections and the potential of
  mechanism-driven modeling,
\newblock \emph{Frontiers in Neuroscience}, \textbf{14} (2020), 1173.

\bibitem{harris2020}
\newblock A.~Harris, G.~Guidoboni, B.~Siesky, S.~Mathew, A.~C.~V. Vercellin,
  L.~Rowe and J.~Arciero,
\newblock Ocular blood flow as a clinical observation: Value, limitations and
  data analysis,
\newblock \emph{Progress in retinal and eye research}, 100841.

\bibitem{harris1996acute}
\newblock A.~Harris, K.~Joos, M.~Kay, D.~Evans, R.~Shetty, W.~E. Sponsel and
  B.~Martin,
\newblock Acute iop elevation with scleral suction: effects on retrobulbar
  haemodynamics.,
\newblock \emph{British journal of ophthalmology}, \textbf{80} (1996),
  1055--1059.

\bibitem{hurtado2017}
\newblock D.~E. Hurtado, S.~Castro and P.~Madrid,
\newblock Uncertainty quantification of 2 models of cardiac electromechanics,
\newblock \emph{International journal for numerical methods in biomedical
  engineering}, \textbf{33} (2017), e2894.

\bibitem{johnson1994lognormal}
\newblock N.~L. Johnson, S.~Kotz and N.~Balakrishnan,
\newblock Lognormal distributions,
\newblock \emph{Continuous univariate distributions}, \textbf{1} (1994),
  601--606.

\bibitem{sala2019cmbe}
\newblock G.~G. L.~Sala~C.~Prud'homme and M.~Szopos,
\newblock The ocular mathematical virtual simulator: towards uncertainty
  quantification,
\newblock in \emph{{\normalfont 6th International Conference on Computational
  and Mathematical Biomedical Engineering (CMBE 2019)} Proceedings}, 2019,
\newblock 429--432.

\bibitem{lee2007association}
\newblock K.~E. Lee, B.~E.~K. Klein, R.~Klein and S.~M. Meuer,
\newblock Association of retinal vessel caliber to optic disc and cup
  diameters,
\newblock \emph{Investigative ophthalmology \& visual science}, \textbf{48}
  (2007), 63--67.

\bibitem{Leguy2011}
\newblock C.~Leguy, E.~Bosboom, A.~Belloum, A.~Hoeks and F.~Van De~Vosse,
\newblock Global sensitivity analysis of a wave propagation model for arm
  arteries,
\newblock \emph{Medical engineering \& physics}, \textbf{33} (2011),
  1008--1016.

\bibitem{limpert2001log}
\newblock E.~Limpert, W.~A. Stahel and M.~Abbt,
\newblock Log-normal distributions across the sciences: keys and clues: on the
  charms of statistics, and how mechanical models resembling gambling machines
  offer a link to a handy way to characterize log-normal distributions, which
  can provide deeper insight into variability and probability—normal or
  log-normal: that is the question,
\newblock \emph{BioScience}, \textbf{51} (2001), 341--352.

\bibitem{petzold1982description}
\newblock L.~R. Petzold,
\newblock \emph{Description of DASSL: a differential/algebraic system solver},
\newblock Technical report, Sandia National Labs., Livermore, CA (USA), 1982.

\bibitem{prieur2017variance}
\newblock C.~Prieur and S.~Tarantola,
\newblock Variance-based sensitivity analysis: Theory and estimation
  algorithms,
\newblock \emph{Handbook of Uncertainty Quantification}, 1217--1239.

\bibitem{quaglino2018}
\newblock A.~Quaglino, S.~Pezzuto, P.-S. Koutsourelakis, A.~Auricchio and
  R.~Krause,
\newblock Fast uncertainty quantification of activation sequences in
  patient-specific cardiac electrophysiology meeting clinical time constraints,
\newblock \emph{International journal for numerical methods in biomedical
  engineering}, \textbf{34} (2018), e2985.

\bibitem{ren2010cerebrospinal}
\newblock R.~Ren, J.~B. Jonas, G.~Tian, Y.~Zhen, K.~Ma, S.~Li, H.~Wang, B.~Li,
  X.~Zhang and N.~Wang,
\newblock Cerebrospinal fluid pressure in glaucoma: a prospective study,
\newblock \emph{Ophthalmology}, \textbf{117} (2010), 259--266.

\bibitem{adams2009dakota}
\newblock \review{Adams, Brian M and Bohnhoff, William J and Dalbey, Keith R
  and Eddy, JP and Eldred, MS and Gay, DM and Haskell, K and Hough, Patricia D
  and Swiler, Laura P},
\newblock \review{DAKOTA, a multilevel parallel object-oriented framework for
  design optimization, parameter estimation, uncertainty quantification, and
  sensitivity analysis: version 5.0 user's manual},
\newblock \emph{\review{Sandia National Laboratories, Tech. Rep.
  SAND2010-2183}}.

\bibitem{Daversin2015}
\newblock \review{Daversin, C. and Prud'homme, C.},
\newblock \review{Simultaneous empirical interpolation and reduced basis method
  for non-linear problems},
\newblock \emph{\review{Comptes Rendus Mathematique}}, \textbf{\review{353}}.

\bibitem{fleeter2020multilevel}
\newblock \review{Fleeter, Casey M. and Geraci, Gianluca and Schiavazzi,
  Daniele E. and Kahn, Andrew M. and Marsden, Alison L.},
\newblock \review{Multilevel and multifidelity uncertainty quantification for
  cardiovascular hemodynamics},
\newblock \emph{\review{Computer Methods in Applied Mechanics and
  Engineering}}, \textbf{\review{365}} (\review{2020}), \review{113030}.

\bibitem{hild:tel-03025312}
\newblock \review{Hild, Romain},
\newblock \emph{\review{Optimization and control of high fields magnets}},
\newblock \review{Theses}, \review{Universit{\'e} de Strasbourg},
  \review{2020},
\newblock \urlprefix\url{https://tel.archives-ouvertes.fr/tel-03025312}.

\bibitem{hose2019cardiovascular}
\newblock \review{Hose, D Rodney and Lawford, Patricia V and Huberts, Wouter
  and Hellevik, Leif Rune and Omholt, Stig W and van de Vosse, Frans N},
\newblock \review{Cardiovascular models for personalised medicine: Where now
  and where next?},
\newblock \emph{\review{Medical engineering \& physics}}, \textbf{\review{72}}
  (\review{2019}), \review{38--48}.

\bibitem{marquis2018practical}
\newblock \review{Marquis, Andrew D and Arnold, Andrea and Dean-Bernhoft, Caron
  and Carlson, Brian E and Olufsen, Mette S},
\newblock \review{Practical identifiability and uncertainty quantification of a
  pulsatile cardiovascular model},
\newblock \emph{\review{Mathematical biosciences}}, \textbf{\review{304}}
  (\review{2018}), \review{9--24}.

\bibitem{Prudhomme2012}
\newblock \review{Prud'homme, C. and Chabannes, V. and Doyeux, V. and Ismail,
  M. and Samake, A. and Pena, G. and Daversin, C. and Trophime, C.},
\newblock \review{Advances in FEEL++ : A domain specific embedded language in
  C++ for partial differential equations},
\newblock in \emph{\review{ECCOMAS 2012 - European Congress on Computational
  Methods in Applied Sciences and Engineering, e-Book Full Papers}},
  \review{2012}.

\bibitem{Prudhomme2020}
\newblock \review{Prud'homme, Christophe and Chabannes, Vincent and Metivet,
  Thibaut and Daversin-Catty, C{\'{e}}cile and Hild, Romain and Doll{\'{e}},
  Guillaume and Sala, Lorenzo and Trophime, Christophe and Samake, Abdoulaye},
\newblock \review{feelpp/feelpp: Feel++ V108},
\newblock \urlprefix\url{https://doi.org/10.5281/zenodo.3784254}.

\bibitem{riva1979bidirectional}
\newblock C.~E. Riva, G.~T. Feke, B.~Eberli and V.~Benary,
\newblock Bidirectional ldv system for absolute measurement of blood speed in
  retinal vessels,
\newblock \emph{Applied Optics}, \textbf{18} (1979), 2301--2306.

\bibitem{riva1985blood}
\newblock C.~E. Riva, J.~E. Grunwald, S.~H. Sinclair and B.~Petrig,
\newblock Blood velocity and volumetric flow rate in human retinal vessels.,
\newblock \emph{Investigative ophthalmology \& visual science}, \textbf{26}
  (1985), 1124--1132.

\bibitem{SundnesCMBE2019}
\newblock S.~T.~W. Rocio Rodriguez-Cantano~Henrik N.~Finsberg and J.~Sundnes,
\newblock A bayesian approach for parameter estimation in computational models
  of cardiac mechanics,
\newblock in \emph{6th International Conference on Computational and
  Mathematical Biomedical Engineering (CMBE 2019) Proceedings P. Nithiarasu, M.
  Ohta, M. Oshima (Eds.)}, 2019,
\newblock 535--538.

\bibitem{sacco2015}
\newblock R.~Sacco, S.~Cassani, G.~Guidoboni, M.~Szopos, C.~Prud'homme and
  A.~Harris,
\newblock Modeling the coupled dynamics of ocular blood flow and production and
  drainage of aqueous humor,
\newblock in \emph{{\normalfont 4th International Conference on Computational
  and Mathematical Biomedical Engineering (CMBE 2015)} Proceedings}, 2015,
\newblock 608--611.

\bibitem{sala2019mathematical}
\newblock L.~Sala,
\newblock \emph{Mathematical modelling and simulation of ocular blood flows and
  their interactions.},
\newblock PhD thesis, Université de Strasbourg, 2019.

\bibitem{sala2018ocular}
\newblock L.~Sala, C.~Prud'Homme, G.~Guidoboni and M.~Szopos,
\newblock Ocular mathematical virtual simulator: A hemodynamical and
  biomechanical study towards clinical applications,
\newblock \emph{Journal of Coupled Systems and Multiscale Dynamics}, \textbf{6}
  (2018), 241--247.

\bibitem{sala2018analysis}
\newblock L.~Sala, C.~Prud'homme, G.~Guidoboni, M.~Szopos, B.~A. Siesky and
  A.~Harris,
\newblock Analysis of iop and csf alterations on ocular biomechanics and lamina
  cribrosa hemodynamics,
\newblock \emph{Investigative Ophthalmology \& Visual Science}, \textbf{59}
  (2018), 4475--4475.

\bibitem{sala2017patient}
\newblock L.~Sala, C.~Prud'Homme, D.~Prada, F.~Salerni, C.~Trophime,
  V.~Chabannes, M.~Szopos, R.~Repetto, S.~Bertoluzza, R.~Sacco et~al.,
\newblock Patient-specific virtual simulator of tissue perfusion in the lamina
  cribrosa.,
\newblock vol.~58,
\newblock The Association for Research in Vision and Ophthalmology, 2017,
\newblock 727.

\bibitem{saltelli2002}
\newblock A.~Saltelli,
\newblock Making best use of model evaluations to compute sensitivity indices,
\newblock \emph{Computer physics communications}, \textbf{145} (2002),
  280--297.

\bibitem{saltelli1999}
\newblock A.~Saltelli, S.~Tarantola and K.-S. Chan,
\newblock A quantitative model-independent method for global sensitivity
  analysis of model output,
\newblock \emph{Technometrics}, \textbf{41} (1999), 39--56.

\bibitem{sankaran2011}
\newblock S.~Sankaran and A.~L. Marsden,
\newblock A stochastic collocation method for uncertainty quantification and
  propagation in cardiovascular simulations,
\newblock \emph{Journal of biomechanical engineering}, \textbf{133}.

\bibitem{sesso2000systolic}
\newblock H.~D. Sesso, M.~J. Stampfer, B.~Rosner, C.~H. Hennekens, J.~M.
  Gaziano, J.~E. Manson and R.~J. Glynn,
\newblock Systolic and diastolic blood pressure, pulse pressure, and mean
  arterial pressure as predictors of cardiovascular disease risk in men,
\newblock \emph{Hypertension}, \textbf{36} (2000), 801--807.

\bibitem{sobol1993sensitivity}
\newblock I.~M. Sobol,
\newblock Sensitivity analysis for non-linear mathematical models,
\newblock \emph{Mathematical modelling and computational experiment},
  \textbf{1} (1993), 407--414.

\bibitem{starling1896absorption}
\newblock E.~H. Starling,
\newblock On the absorption of fluids from the connective tissue spaces,
\newblock \emph{The Journal of physiology}, \textbf{19} (1896), 312--326.

\bibitem{suh2012distribution}
\newblock W.~Suh, C.~Kee, N.~S. Group and K.~G. Society,
\newblock The distribution of intraocular pressure in urban and in rural
  populations: the namil study in south korea,
\newblock \emph{American journal of ophthalmology}, \textbf{154} (2012),
  99--106.

\bibitem{szopos2016}
\newblock M.~Szopos, S.~Cassani, G.~Guidoboni, C.~Prud'homme, R.~Sacco,
  B.~Siesky and A.~Harris,
\newblock {Mathematical modeling of aqueous humor flow and intraocular pressure
  under uncertainty: towards individualized glaucoma management},
\newblock \emph{J. for Modeling in Ophthalmology}, \textbf{1} (2016), 29--39.

\bibitem{tham2018inter}
\newblock Y.-C. Tham, S.-H. Lim, P.~Gupta, T.~Aung, T.~Y. Wong and C.-Y. Cheng,
\newblock Inter-relationship between ocular perfusion pressure, blood pressure,
  intraocular pressure profiles and primary open-angle glaucoma: the singapore
  epidemiology of eye diseases study,
\newblock \emph{British Journal of Ophthalmology}, \textbf{102} (2018),
  1402--1406.

\bibitem{vercellin2016}
\newblock A.~C.~V. Vercellin, A.~Harris, J.~V. Cordell, T.~Do, J.~Moroney,
  A.~Belamkar and B.~Siesky,
\newblock Mathematical modeling and glaucoma: the need for an individualized
  approach to risk assessment,
\newblock \emph{Journal for Modeling in Ophthalmology}, \textbf{1} (2016),
  6--20.

\bibitem{wang2011intraocular}
\newblock D.~Wang, W.~Huang, Y.~Li, Y.~Zheng, P.~J. Foster, N.~Congdon and
  M.~He,
\newblock Intraocular pressure, central corneal thickness, and glaucoma in
  chinese adults: the liwan eye study,
\newblock \emph{American journal of ophthalmology}, \textbf{152} (2011),
  454--462.

\bibitem{xu2011reliability}
\newblock C.~Xu and G.~Z. Gertner,
\newblock Reliability of global sensitivity indices,
\newblock \emph{Journal of Statistical Computation and Simulation}, \textbf{81}
  (2011), 1939--1969.

\end{thebibliography}


\end{document}